\documentclass[notitlepage,10pt]{article}
\usepackage{amssymb,mathrsfs,amsmath,enumitem}
\catcode`\@=11
\@addtoreset{equation}{section}

\catcode`\@=12

\newtheorem{Theorem}{Theorem}[section]
\newtheorem{Lemma}[Theorem]{Lemma}

\newtheorem{Remark}[Theorem]{Remark}
\newtheorem{Definition}[Theorem]{Definition}

\def\R{{\mathbb R}}

\def\N{{\mathbb N}}
\def\S{{\mathbb S}}

\def\A{{\mathcal A}}

\def\D{{\mathcal D}}
\def\E{{\mathcal E}}
\def\F{{\mathcal F}}
\def\G{{\mathcal G}}
\def\Q{{\mathcal Q}}
\def\V{{\mathcal V}}
\def\disc{\mathbb{D}}
\def\Proof{\noindent\textit{Proof. }}
\def\qed{$~\square$\goodbreak \medskip}

\title{Isovolumetric and isoperimetric problems\\
for a class of capillarity functionals}
\author{Paolo Caldiroli\footnote{Dipartimento di Matematica, Universit\`a di Torino, via Carlo Alberto, 10 -- 10123 Torino, Italy. Email: \tt{paolo.caldiroli@unito.it}}}

\begin{document}
%\centerline{\Large{\textbf{Isovolumetric and isoperimetric problems}}}
%\medskip
%
%\centerline{\Large{\textbf{for a class of capillarity functionals}}}
\date{}
\maketitle
%\vspace{-50pt}

\begin{abstract}
\noindent
Capillarity functionals are parameter invariant functionals defined on classes of two-dimensional parametric surfaces in $\R^{3}$ as the sum of the area integral and an anisotropic term of suitable form. In the class of parametric surfaces with the topological type of $\S^{2}$ and with fixed volume, extremals of capillarity functionals are surfaces whose mean curvature is prescribed up to a constant. For a certain class of anisotropies vanishing at infinity, we prove existence and nonexistence of volume-constrained, $\S^{2}$-type, minimal surfaces for the corresponding capillarity functionals. Moreover, in some cases, we show existence of extremals for the full isoperimetric inequality. 
\smallskip

\noindent
\textit{Keywords:} {Isovolumetric problems, isoperimetric problems, parametric surfaces, $H$-bubbles.}
\smallskip

\noindent{{{\it 2010 Mathematics Subject Classification:} 53A10 (49Q05, 49J10)}}
\end{abstract}

\section{Introduction}

In this work we deal with closed surfaces in $\R^{3}$ parametrized by mappings $u\colon\S^{2}\to\R^{3}$. Introducing the stereographic projection $\phi$ of $\S^{2}$ onto the compactified plane $\R^{2}\cup\{\infty\}$ and identifying maps $u$ defined on $\S^{2}$ with corresponding maps $u\circ\phi^{-1}$ on $\R^{2}\cup\{\infty\}$, the area of a surface parametrized by $u$ is given by
$$
\A(u):=\int_{\R^{2}}|u_{x}\wedge u_{y}|
$$
whereas the algebraic volume enclosed by $u$ can be computed in terms of the Bononcini-Wente integral
$$
\V(u):=\frac{1}{3}\int_{\R^{2}}u\cdot u_{x}\wedge u_{y}.
$$
The relationship between the area and the volume integrals is stated by the classical isoperimetric inequality, proved in \cite{Bon53}:
\begin{equation}
\label{eq:ii}
S|\V(u)|^{2/3}\le\A(u)\quad\forall u\in C^{\infty}(\S^{2},\R^{3})
\end{equation}
where $S=\sqrt[3]{36\pi}$. As one expects, the inequality (\ref{eq:ii}) in fact holds true in the Sobolev space $H^{1}(\S^{2},\R^{3})$ (see \cite{Wen69}) and the constant $S=\sqrt[3]{36\pi}$ is the best one and is achieved if and only if $u$ parametrizes a round sphere with arbitrary center and radius. (This fact can be readily deduced from the results discussed in \cite{BreCor85}, in particular Lemma 0.1. For a self-contained and direct proof, see \cite{CaMuARMA}, Lemma 2.1.) 

The area integral $\A(u)$ constitutes the simplest and most relevant example of a Cartan functional. As displayed in \cite{DHS}, Sect.~\!4.13, these are integrals of the kind
$$
\F(u):=\int_{\R^{2}}F(u,u_{x}\wedge u_{y})
$$
with a Lagrangian $F\in C^{0}(\R^{3}\times\R^{3})$ such that:
\begin{itemize}[leftmargin=24pt]
\item[$(C_{1})$] $F(p,q)$ is positively homogeneous of degree one with respect to $q$, i.e., $F(p,tq)=tF(p,q)$ for $t>0$ and for all $(p,q)\in\R^{3}\times\R^{3}$,
\item[$(C_{2})$] there exist $0<m_{1}\le m_{2}$ such that the definiteness condition $m_{1}|q|\le F(p,q)\le m_{2}|q|$ holds for all  $(p,q)\in\R^{3}\times\R^{3}$,
\item[$(C_{3})$]
$F(p,q)$ is weakly elliptic, namely it is convex with respect to $q$, i.e. $F(p,tq_{1}+(1-t)q_{2})\le tF(p,q_{1})+(1-t)F(p,q_{2})$ for $t\in[0,1]$ and $p,q_{1},q_{2}\in\R^{3}$. 
\end{itemize}
By $(C_{1})$ and the upper bound in $(C_{2})$ any Cartan functional $\F$ turns out to be well defined in $H^{1}(\S^{2},\R^{3})$ and is a parameter invariant integral, i.e., we have $\F(u\circ g)=\F(u)$ for any $C^{1}$ diffeomorphism $g$ of $\S^{2}$ onto itself. This rightly reflects the geometrical character of the problem we deal with. 

We point out that the framework described above admits a counterpart in the setting of Geometric Measure Theory. In that context, surfaces are meant as boundaries of sets of finite perimeter and Cartan functionals are replaced by boundary functionals defined by so-called ``semielliptic'' integrals (see \cite{CafDel01}, Sect.~\!2). Later we will come back to this aspect.

Thanks to the lower positive bound in $(C_{2})$, and by (\ref{eq:ii}),  an isoperimetric-like inequality for any Cartan functional can be also written, i.e.,
\begin{equation}
\label{eq:ii-Cartan}
S_{F}|\V(u)|^{2/3}\le\F(u)\quad\forall u\in H^{1}(\S^{2},\R^{3})
\end{equation}
for some constant $S_{F}\in(0,m_{1}S]$. The existence of extremals for (\ref{eq:ii-Cartan}) arises as a natural question and constitutes a rather challenging target. Indeed, since in general a Cartan functional is not purely quadratic, differently from (\ref{eq:ii}), the inequality (\ref{eq:ii-Cartan}) is not invariant under dilation and translation (with respect to $u$). These missing invariances might make difficult restoring some compactness for sequences of approximate extremals of (\ref{eq:ii-Cartan}). 

A way to prevent, hopefully, troubles due to dilation is to consider isovolumetric problems, i.e., constrained minimimization problems with fixed volume, as follows. Fixing $t\in\R$, study the existence of minimizers for
\begin{equation}
\label{eq:inf-Cartan}
S_{F}(t):=\inf\{\F(u)~|~u\in H^{1}(\S^{2},\R^{3}),~\V(u)=t\}.
\end{equation}
We point out that also these minimization problems are far from being obvious because, even if Cartan functionals are weakly lower semicontinuous (see \cite{DHS}), the constraint is not weakly closed and the volume functional is not weakly lower semicontinuous (see \cite{Wen69}). 
In fact, as we will see in some cases, the existence or nonexistence of minimizers for (\ref{eq:inf-Cartan}) is a rather delicate issue and depends in a sensitive way on the shape of the Lagrangian.

In this paper we study problems (\ref{eq:inf-Cartan}) for a special class of Lagrangian functions. In particular we consider 
%\begin{equation}\label{eq:capillarity}
$$F(p,q)=|q|+Q(p)\cdot q$$
%\end{equation}
with $Q\in C^{1}(\R^{3},\R^{3})$ prescribed, such that $\|Q\|_{\infty}<1$. Cartan functionals corresponding to such $F$, which indeed satisfy $(C_{1})$--$(C_{3})$, can be interpreted as modified area integrals with an anisotropy term:
%\begin{equation}\label{eq:capillarity-functional}
$$\F(u)=\A(u)+\int_{\R^{2}}Q(u)\cdot u_{x}\wedge u_{y}$$
%\end{equation}
and are often known as ``capillarity functionals'' (see \cite{HilMos}). They are particularly meaningful because in this case possible minimizers for (\ref{eq:inf-Cartan}) parametrize $\S^{2}$-type surfaces with volume $t$ and mean curvature $H(p)=K(p)-\lambda$ where $K=\mathrm{div}~\!Q$ is  prescribed, whereas $\lambda$ is a constant corresponding to the Lagrange multiplier due to the constraint. We will call such surfaces ``$H$-bubbles''. In the sequel the strong relation between the isovolumetric problem for capillarity functionals and the $H$-bubble problem will become even more evident. 

Capillarity functionals depend on the vector field $Q$ only by its divergence. Therefore we can state the precise assumptions just on the scalar field $K=\mathrm{div}~\!Q$. In the present work we focus on a class of mappings $K\colon\R^{3}\to\R$ vanishing at infinity with a suitable rate. In particular let us start by assuming that $K\in C^{1}(\R^{3})$ satisfies: 
\begin{itemize}[leftmargin=25pt]
\item[$(K_{1})$]~~$|K(p)p|\le k_{0}<2$ for every $p\in\R^{3}$.
\item[$(K_{2})$]~~$K(p)p\to 0$ as $|p|\to\infty$. 
\end{itemize}
Then it is possible to construct a vector field $Q_{K}\in C^{1}(\R^{3},\R^{3})$ such that $\mathrm{div}~\!Q_{K}=K$ on $\R^{3}$ and enjoying the following properties:
\begin{itemize}[leftmargin=25pt]
\item[$(Q_{1})$]~~$\|Q_{K}\|_{\infty}<1$,
\item[$(Q_{2})$]~~$|Q_{K}(p)|\to 0$ as $|p|\to\infty$. 
\end{itemize}
These are direct consequences of $(K_{1})$ and $(K_{2})$, respectively  (see Remark \ref{R:K-volume}). Therefore the assumptions $(K_{1})$ and $(K_{2})$ seem to be reasonably natural to deal with situations with anisotropy vanishing at infinity. 

In order to state a satisfactory result about the minimization problems
\begin{equation}
\label{eq:SKt}
\begin{split}
S_{K}(t):=\inf\left\{\F_{K}(u)~|~u\in H^{1}(\S^{2},\R^{3}),~\V(u)=t\right\}\\
\text{where}\quad\F_{K}(u):=\A(u)+\int_{\R^{2}}Q_{K}(u)\cdot u_{x}\wedge u_{y},~~
\end{split}
\end{equation}
in addition to the conditions $(K_{1})$ and $(K_{2})$, we introduce an extra assumption which controls the radial oscillation of $K$:
\begin{itemize}[leftmargin=25pt]
\item[$(K_{3})$]~~$|(\nabla K(p)\cdot p)p|\le k_{0}<2~~\forall p\in\R^{3}$.
\end{itemize}
We point out that $(K_{3})$ together with $(K_{2})$ implies $(K_{1})$ (see \cite{CaMu11}, Remark 2.2, for a proof). The first existence result shown in this paper can be stated as follows. 

\begin{Theorem}
\label{T:existence1}
Let $K\in C^{1}(\R^{3})$ satisfy $(K_{2})$ and $(K_{3})$. Let
\begin{equation}
\label{eq:tpm}
\begin{split}
&t_{+}:=\sup\left\{t\ge 0~|~K\le 0\text{ and }K\not\equiv 0\text{ in some ball of radius }\sqrt[3]{3t/4\pi}\right\}\\
&t_{-}:=\inf\left\{t\le 0~|~K\ge 0\text{ and }K\not\equiv 0\text{ in some ball of radius }\sqrt[3]{3|t|/4\pi}\right\}.
\end{split}
\end{equation}
Then for every $t\in(t_{-},t_{+})$ there exists $U\in H^{1}(\S^{2},\R^{3})$ with $\V(U)=t$ and $\F_{K}(U)=S_{K}(t)$. Moreover when $t\ne 0$ such $U$ is a $(K-\lambda)$-bubble, of class $C^{2,\alpha}$, for some $\lambda=\lambda(t,U)\ne 0$. 
\end{Theorem}
Notice that, in the definition of $t_{+}$ and $t_{-}$, one could have balls with arbitrary (and in general different) centers. Moreover, excluding the trivial case $K=0$, the interval $(t_{-},t_{+})$ is always nonempty. 

In fact, the sign of $K$ plays a crucial role in the above stated result. For example, if $K<0$ (respectively, $K>0$) on the tail of some open cone, then $t_{+}=\infty$ (resp., $t_{-}=-\infty$). 

It is not clear if the result stated in Theorem \ref{T:existence1} is optimal. But in some cases we can provide some more information. In particular, when $K<0$ on $\R^{3}$, then, according to Theorem \ref{T:existence1}, a minimizer for problem (\ref{eq:SKt}) exists for every $t>0$. Actually, we can show non existence of minimizers as $t<0$, but just for small $|t|$ (see Theorem \ref{T:nonexistence}). 

The arguments of the proof make full use of refined tools already developed in the context of the $H$-bubble problem. In particular the study of minimizing sequences for the isovolumetric problems (\ref{eq:SKt}) exploits some deep results proved in \cite{CaMuJFA} and \cite{CaJFA}, concerning the behavior of approximate solutions of $(K-\lambda)$-systems
\begin{equation}
\label{eq:K-lambda-system}
\Delta u=(K(u)-\lambda)u_{x}\wedge u_{y}\quad\text{on}\quad\R^{2}
\end{equation}
which, to our knowledge, are known just when the mapping $K$ satisfies precisely $(K_{1})$ and $(K_{2})$. 

In fact, all the assumptions asked of $K$ in Theorem \ref{T:existence1} are the same considered in the papers \cite{CaMuCCM} and \cite{CaMu11} on the $H$-bubble problem for a prescribed mean curvature function $H(p)=H_{0}(p)+H_{\infty}$ where $H_{\infty}$ is a nonzero constant corresponding to $-\lambda$, whereas $H_{0}$ is a $C^{1}$ function on $\R^{3}$, vanishing at infinity, playing a similar role of $K$. The only difference is a factor $2$, because in \cite{CaMuCCM} and \cite{CaMu11} one writes the prescribed mean curvature equation for parametric surfaces in the form $\Delta u=2H(u)u_{x}\wedge u_{y}$. Conditions $(K_{1})$ and $(K_{3})$ are changed accordingly.

Clearly, for the $H$-bubble problem the volume of the solution is not prescribed. Moreover, in the works \cite{CaMuCCM} and \cite{CaMu11}, solutions are found as saddle-type critical points of the (unbounded) energy functional naturally associated to (\ref{eq:K-lambda-system}). Furthermore, in general, nonconstant weak solutions $u\in H^{1}(\S^{2},\R^{3})$ of (\ref{eq:K-lambda-system}) are not necessarily minimizers for the isovolumetric problem (\ref{eq:SKt}) with $t=\V(u)$. 

Considering the set of mean curvature functions $H=K-\lambda$ for which Theorem \ref{T:existence1} provides existence of a $(K-\lambda)$-bubble, we cannot guarantee that the range of admissible values for $\lambda$ does not contain gaps. On the other hand, the occurrence of gaps would not be surprising, when the metric induced by the anisotropy term is far from flat (see \cite{Ban87} and \cite{GolNov12} for examples in this spirit, but in different contexts). Anyway, some information on the set-valued function 
\begin{equation*}
\begin{split}
t\mapsto\Lambda(t):=\{\lambda\in\R~|~\exists U\in H^{1}(\S^{2},\R^{3})~&\text{minimizer of $S_{K}(t)$}\\
&\text{and $(K-\lambda)$-bubble}\}\quad(t\in(t_{-},t_{+}))
\end{split}
\end{equation*}
is available (see Theorem \ref{T:Lambda}). 

A few more words can be said about the assumption $(K_{3})$. This condition, which is essential in the works \cite{CaMuCCM} and \cite{CaMu11} about the $H$-bubble problem, here plays a role just in order to avoid that minimizing sequences for the isovolumetric problems (\ref{eq:SKt}) split into many $(K-\lambda)$-bubbles (see Lemma \ref{L:strict-inequality}). It is not clear if $(K_{3})$ is a purely technical assumption. As a matter of fact, we can provide a second existence result for the isovolumetric problems without $(K_{3})$, just assuming $(K_{1})$ and $(K_{2})$, but with a restriction on the constant $k_{0}$ appearing in $(K_{1})$. More precisely, we have:

\begin{Theorem}
\label{T:existence2}
Let $K\in C^{1}(\R^{3})$ satisfy $(K_{1})$, $(K_{2})$, and 
\begin{itemize}[leftmargin=25pt]
\item[$(K_{4})$]~$2^{2/3}(2+k_{0})<(2-k_{0})^{2}$.
\end{itemize}
Moreover, let $t_{+}$ and $t_{-}$ be defined as in (\ref{eq:tpm}). 
Then the same conclusion of Theorem \ref{T:existence1} holds true.
\end{Theorem}

Condition $(K_{4})$, even if somehow unnatural, is worth considering because it does not involve derivatives of $K$. Furthermore, thanks to Theorem \ref{T:existence2} and to the information about Lagrange multipliers $\lambda=\lambda(t)$, we obtain a new result about existence of $H$-bubbles with prescribed mean curvature $H$ assuming a large constant value at infinity (see Theorem \ref{T:existence4}). 

In the second part of this work we turn attention to isoperimetric inequalities like (\ref{eq:ii-Cartan}) for capillarity functionals $\F_{K}$ with $K$ of the form considered before and, pushing on the investigation, we prove the following existence result. 

\begin{Theorem}
\label{T:existence3}
Let $K\in C^{1}(\R^{3})$ satisfy $(K_{2})$--$(K_{3})$ or, as an alternative, $(K_{1})$, $(K_{2})$, and $(K_{4})$. If $K\le 0$ on $\R^{3}$ then, letting $\F_{K}$ as in (\ref{eq:SKt}), the minimization problem
\begin{equation}
\label{eq:K-isoperimetric-inequality}
S_{K}:=\inf_{\scriptstyle{u\in H^{1}(\S^{2},\R^{3})}\atop\scriptstyle{\V(u)>0}}\frac{\F_{K}(u)}{\V(u)^{2/3}}
\end{equation}
admits a solution. Moreover if $U\in H^{1}(\S^{2},\R^{3})$ is a minimizer for (\ref{eq:K-isoperimetric-inequality}), then $U$ is a $(K-\lambda)$-bubble, of class $C^{2,\alpha}$, with $\lambda=\frac{2}{3}S_{K}\V(U)^{-1/3}$.
\end{Theorem}

As mentioned at the beginning, isovolumetric-type problems, like those considered in this paper, might be tackled also using methods of Geometric Measure Theory. For example, this is carried out in \cite{CafDel01} and \cite{GolNov12} in case of periodic media. 

However, we would like to stress that we are interested in volume-con\-strained minimal surfaces with the topological type of the sphere. A geometric measure-theoretic approach seems to lack in providing this kind of information whereas, under global assumptions on the anisotropy, the approach by means of pa\-ram\-e\-tri\-za\-tions, as followed in this work, turns out to be well suited to this purpose. 

Moreover, we expect that the general structure displayed here could be hopefully adapted in dealing with different, maybe more general, classes of anisotropies and, in a wider perspective, could be possibly constitute an alternative method to tackle the $H$-bubble problem.

\section{Preliminaries}
\label{S:preliminaries}
Let us introduce the space
$$
\hat{H}^{1}:=\{u\in H^{1}_{loc}(\R^{2},\R^{3})~|~\int_{\R^{2}}(|\nabla u|^{2}+\mu^{2}|u|^{2})<\infty\}
$$ 
where 
\begin{equation}
\label{eq:mu}
\mu(z)=\frac{2}{1+|z|^{2}}\quad\text{for } z\in\R^{2}. 
\end{equation}
The space $\hat{H}^{1}$ is a Hilbert space with inner product
$$
\langle u,v\rangle=\int_{\R^{2}}(\nabla u\cdot\nabla v)+\left(\frac{1}{4\pi}\int_{\R^{2}}u\mu^{2}\right)\cdot\left(\frac{1}{4\pi}\int_{\R^{2}}v\mu^{2}\right)
$$
and is isomorphic to the space $H^{1}(\S^{2},\R^{3})$. The isomorphism is given by the correspondence $\hat{H}^{1}\ni u\mapsto u\circ\phi\in H^{1}(\S^{2},\R^{3})$, where $\phi$ is the stereographic projection of $\S^{2}$ onto the compactified plane $\R^{2}\cup\{\infty\}$. As usual, we denote $\|u\|=\langle u,u\rangle^{1/2}$.

One has that $C^{\infty}(\S^{2},\R^{3})$ is dense in $H^{1}(\S^{2},\R^{3})$ (see, e.g., \cite{Aub}, Ch.2). As a consequence, 
$\hat{C}^{\infty}:=\{u\circ\phi^{-1}~|~u\in C^{\infty}(\S^{2},\R^{3})\}$ is dense in $\hat{H}^{1}$. We point out that $p+\hat{H}^{1}=\hat{H}^{1}$ for every $p\in\R^{3}$. By obvious extension, for every bounded domain $\Omega$ in $\R^{2}$ the space $H^{1}_{0}(\Omega,\R^{3})$ can be considered as a subspace of $\hat{H}^{1}$. Then also $p+H^{1}_{0}(\Omega,\R^{3})$ is an affine subspace of $\hat{H}^{1}$ for every $p\in\R^{3}$. 

\begin{Lemma}
\label{L:density}
The space $\R^{3}+C^{\infty}_{c}(\R^{2},\R^{3})$ is dense in $\hat{H}^{1}$. In particular, for every $u\in\hat{H}^{1}\cap L^{\infty}$ there exists a sequence $(u^{n})\subset\R^{3}+C^{\infty}_{c}(\R^{2},\R^{3})$ such that $u^{n}\to u$ in $\hat{H}^{1}$, in $L^{\infty}_{loc}$ and $\|u^{n}\|_{\infty}\le\|u\|_{\infty}$.
\end{Lemma}

Even if this result is known, for future convenience, we sketch a proof, which contains a construction used also in the sequel.
\medskip

\noindent
\Proof
Take $u\in\hat{C}^{\infty}$, let $p=\lim_{|z|\to\infty}u(z)$, and for every $n\in\N$ set
\begin{equation}
\label{eq:approx-u}
u^{n}(z)=\left\{\begin{array}{ll}u(z)&\text{as }|z|\le n\\ \left(2-\frac{\log|z|}{\log n}\right)u(z)+\left(\frac{\log|z|}{\log n}-1\right)p&\text{as }n<|z|\le n^{2}\\ p&\text{as }|z|>n^{2}.\end{array}\right. 
\end{equation}
Setting $A_{n}=\{z\in\R^{2}~|~n<|z|\le n^{2}\}$, we have that
\begin{equation*}
\begin{split}
\int_{\R^{2}}|\nabla(u-u^{n})|^{2}&=\int_{A_{n}}\left|\nabla\left[\left(1-\frac{\log|z|}{\log n}\right)\left(u(z)-p\right)\right]\right|^{2}+\int_{|z|>n^{2}}|\nabla u|^{2}\\
&\le\int_{|z|>n}|\nabla u|^{2}+\sup_{|z|>n}|u(z)-p|\int_{A_{n}}\left|\frac{\nabla(\log|z|)}{\log n}\right|^{2}=o(1)
\end{split}
\end{equation*}
as $n\to\infty$. Moreover
$$
\left|\int_{\R^{2}}(u-u^{n})\mu^{2}\right|\le\sup_{|z|>n}|u(z)-p|\int_{|z|>n}\mu^{2}=o(1)
$$
as $n\to\infty$. For every $\varepsilon>0$ there exists $n\in\N$ such that $\|u-u^{n}\|<\frac{\varepsilon}{2}$. Since $u^{n}-p\in H^{1}_{0}(\Omega_{n},\R^{3})$, where $\Omega_{n}$ is the disc of radius $n^{2}$, we can find $v\in C^{\infty}_{c}(\Omega_{n})$ such that $\|u^{n}-v\|<\frac{\varepsilon}{2}$. Hence the conclusion follows from the density of $\hat{C}^{\infty}$ in $\hat{H}^{1}$. Since this last property can be proved by a standard regularizing technique using Friedrichs mollifiers which do not increase the $L^{\infty}$ norm, also the second part of the lemma is true. 
\qed\medskip

Set
$$
\D(u):=\frac{1}{2}\int_{\R^{2}}|\nabla u|^{2}\quad(u\in\hat{H}^{1})~~\text{and}~~\V(u):=\frac{1}{3}\int_{\R^{2}}u\cdot u_{x}\wedge u_{y}\quad(u\in\hat{H}^{1}\cap L^{\infty}).
$$

\begin{Lemma}
\label{L:volume}
The functional $\V$ admits a unique analytic extension on $\hat{H}^{1}$. In particular for every $u\in\hat{H}^{1}$
$$
\V'(u)[\varphi]=\int_{\R^{2}}\varphi\cdot u_{x}\wedge u_{y}\quad\forall\varphi\in \hat{H}^{1}\cap L^{\infty}
$$
and there exists a unique $v\in\hat{H}^{1}\cap L^{\infty}$ which is a (weak) solution of 
\begin{equation}
\label{eq:V-prime-representative}
\left\{\begin{array}{l}
-\Delta v=u_{x}\wedge u_{y}\quad\text{on }\R^{2}\\
\int_{\R^{2}}v\mu^{2}=0.
\end{array}\right.
%,\quad\V'(u)=\langle v,\cdot\rangle.
\end{equation}
Moreover 
\begin{equation}
\label{eq:Wente-inequality}
\|\nabla v\|_{2}+\|v\|_{\infty}\le C\|\nabla u\|_{2}^{2}
\end{equation}
for a constant $C$ independent of $u$. 
In addition, for every $t\ne 0$ the set 
\begin{equation}
\label{eq:Mt-def}
M_{t}:=\{u\in\hat{H}^{1}~|~\V(u)=t\}
\end{equation}
is a smooth manifold and, for any fixed $u\in M_{t}$, a function $\varphi\in\hat{H}^{1}$ belongs to the tangent space to $M_{t}$ at $u$, denoted $T_{u}M_{t}$, if and only if $\V'(u)[\varphi]=0$. 
%is $T_{u}M_{t}=\{\varphi\in\hat{H}^{1}~|~\V'(u)[\varphi]=0\}$. 
\end{Lemma}

\begin{Remark}
\label{R:volume}
The second part of Lemma \ref{L:volume} states that there exists $C>0$ such that $\|\V'(u)\|_{\hat{H}^{-1}}\le C\|\nabla u\|_{2}^{2}$ for every $u\in\hat{H}^{1}$, where $\hat{H}^{-1}$ denotes the dual of $\hat{H}^{1}$. 
\end{Remark}

\Proof
All the results stated in the lemma are essentially well known; the proof displayed, e.g., in \cite{Wen69}, Thms.~\!3.1 and 3.3 (see also \cite{Str89-book}, Ch.~\!III, Thm.~\!2.3), considering as a domain the space $H^{1}_{0}(\disc,\R^{3})$, where $\disc$ denotes the unit disc in $\R^{2}$, works also in $\hat{H}^{1}$. The only additional remark regards the fact that, for fixed $u\in\hat{H}^{1}$, the Riesz representative of $\V'(u)$ in $\hat{H}^{1}$ belongs to $L^{\infty}$. To prove this, we consider a sequence of Dirichlet problems
\begin{equation}
\label{eq:Wente-n}
\left\{\begin{array}{ll}-\Delta v=u_{x}\wedge u_{y}&\text{in }\Omega_{n}\\ v=0&\text{on }\partial\Omega_{n}\end{array}\right.
\end{equation}
where $\Omega_{n}=\{z\in\R^{2}~|~|z|<n\}$. It is known that for every $n\in\N$ there exists $v^{n}\in H^{1}_{0}$ solving (\ref{eq:Wente-n}) and 
$$
\|\nabla v^{n}\|_{2}+\|v^{n}\|_{\infty}\le C\|\nabla u\|^{2}_{2}
$$
with $C$ independent of $n$ (see \cite{BetGhi93}; see also \cite{Top97} for the optimal constant $C=(2\pi)^{-1}$). Then the sequence $(v^{n})$ is bounded in $\hat{H}^{1}$ and in $L^{\infty}$, admits a subsequence which converges weakly in $\hat{H}^{1}$ to some $w\in\hat{H}^{1}\cap L^{\infty}$ solving
$$
\int_{\R^{2}}\nabla w\cdot\nabla\varphi=\int_{\R^{2}}\varphi\cdot u_{x}\wedge u_{y}\quad\forall\varphi\in \R^{3}+C^{\infty}_{c}(\R^{2},\R^{3})
$$
(notice that $\int_{\R^{2}}p\cdot u_{x}\wedge u_{y}=0$ for all $p\in\R^{3}$). Finally, the function $v=w-\frac{1}{4\pi}\int_{\R^{2}}w\mu^{2}$ solves (\ref{eq:V-prime-representative}) and belongs to $L^{\infty}$. 
\qed

\begin{Remark}
\label{R:sphere}
The mapping $\omega(z)=\left(\mu x,\mu y,1-\mu\right)$, with $\mu$ defined in (\ref{eq:mu}), is a conformal parame\-tri\-za\-tion of the unit sphere. Indeed, it is the inverse of the stereographic projection from the North Pole.  Moreover $\A(\omega)=\D(\omega)=4\pi$ and $\V(\omega)=-\frac{4\pi}{3}$. If $\overline{p}\in\R^{3}$ and $r\in\R\setminus\{0\}$, then $u=\overline{p}+r\omega$ is a parametrization of a sphere centered at $\overline{p}$ and with radius $|r|$, Moreover $\A(u)=\D(u)=4\pi r^{2}$ and $\V(u)=-\frac{4\pi r^{3}}{3}$.
\end{Remark}

\begin{Lemma}[Isoperimetric inequality]
\label{L:isoperimetric}
It holds that 
\begin{equation}
\label{eq:isoperimetric}
S|\V(u)|^{2/3}\le\A(u)\le\D(u)\quad\forall u\in \hat{H}^{1}
\end{equation}
where $S=\sqrt[3]{36\pi}$ is the best constant. Moreover any extremal function for (\ref{eq:isoperimetric}) is a conformal parametrization of a simple sphere. 
%of degree $1$ or $-1$ (as a map from $\S^{2}$ into $\S^{2}$). 
\end{Lemma}
Inequality (\ref{eq:isoperimetric}) for regular mappings goes back to \cite{Bon53}. Its extension to $H^{1}_{0}(\disc,\R^{3})$ is proved in \cite{Wen69}. The version for mappings in $\hat{H}^{1}$, even not explicitly stated, can be also deduced from \cite{Wen69}, Theorem 2.5, by a density argument, by means of Lemmas \ref{L:density} and \ref{L:volume}.  
\medskip

Fixing $K\in C^{1}(\R^{3})$ satisfying $(K_{1})$, set 
$$
m_{K}(p):=\int^{1}_{0}K(sp)s^{2}~\!ds\quad\text{and}\quad Q_{K}(p):=m_{K}(p)p\quad\forall p\in\R^{3}
$$
and observe that $\mathrm{div}~\!Q_{K}=K$. 
Then set
$$
\Q(u):=\int_{\R^{2}}Q_{K}(u)\cdot u_{x}\wedge u_{y}\quad(u\in\hat{H}^{1}).
$$

\begin{Remark}
\label{R:K-volume}
From $|K(p)p|\le k_{0}$ for every $p\in\R^{3}$, it follows that $\|Q_{K}\|_{\infty}\le\frac{k_{0}}{2}$. In particular, by the assumption $(K_{1})$, 
\begin{equation}
\label{eq:Q<1}
\|Q_{K}\|_{\infty}<1.
\end{equation}
Moreover the functional $\Q$ is well defined on $\hat{H}^{1}$ and 
\begin{equation}
\label{eq:VK-estimate}
|\Q(u)|\le\|Q_{K}\|_{\infty}\D(u)\quad\forall u\in\hat{H}^{1}.
\end{equation}
One can also check that 
\begin{equation}
\label{eq:Q(infty)=0}
|Q_{K}(p)|\to 0\quad\text{as}\quad |p|\to\infty.
\end{equation}
Indeed, for $|p|>R$ write 
$$
Q_{K}(p)=\frac{\hat{p}}{|p|^{2}}\int_{0}^{R}K(t\hat{p})t^{2}dt+\frac{\hat{p}}{|p|^{2}}\int_{R}^{|p|}K(t\hat{p})t^{2}dt
$$
with $\hat{p}=\frac{p}{|p|}$, and use $(K_{2})$ to conclude. 
\end{Remark}

The next result collects some useful properties of the functional $\Q$.

\begin{Lemma}
\label{L:K-volume}
Let $K\colon\R^{3}\to\R$ be a bounded continuous function. Then:
\begin{itemize}[leftmargin=22pt]
\item[(i)]
the functional $\Q$ is continuous in $\hat{H}^{1}$.
\item[(ii)]
For every $u\in\hat{H}^{1}$ and $\varphi\in\hat{H}^{1}\cap L^{\infty}$ one has 
$$
\Q(u+\varphi)-\Q(u)=\int_{\R^{2}}\int_{0}^{1}K(u+r\varphi)\varphi\cdot(u_{x}+r\varphi_{x})\wedge(u_{y}+r\varphi_{y})~\!dr~\!dz.
$$
\item[(iii)]
The functional $\Q$ admits directional derivatives at every $u\in\hat{H}^{1}$ along any $\varphi\in\hat{H}^{1}\cap L^{\infty}$, given by
$$
\Q'(u)[\varphi]=\int_{\R^{2}}K(u)\varphi\cdot u_{x}\wedge u_{y}.
$$
\end{itemize}
If in addition $\sup_{p\in\R^{3}}|K(p)p|<\infty$ then for every $u\in\hat{H}^{1}$ the mapping $s\mapsto\Q(su)$ is differentiable and
\begin{equation}
\label{eq:radial-VK-derivative}
\frac{d}{ds}[\Q(su)]=s^{2}\int_{\R^{2}}K(su)u\cdot u_{x}\wedge u_{y}.
\end{equation}
\end{Lemma}

\Proof
The first part of the lemma is proved in \cite{Ste76}, Proposition 3.3, whereas (\ref{eq:radial-VK-derivative}) is discussed in \cite{CaMu11} (in particular, see formula (2.7) therein). 
\qed

\begin{Remark}
\label{R:sphere2}
Let $\omega$ be the mapping introduced in Remark \ref{R:sphere}. For every $\overline{p}\in\R^{3}$ and $r>0$, one has that $\Q(\overline{p}+r\omega)=-\int_{B_{r}(\overline{p})}K(p)~\!dp$ whereas if $r<0$ then $\Q(\overline{p}+r\omega)=\int_{B_{|r|}(\overline{p})}K(p)~\!dp$.
\end{Remark}

We conclude this section with an auxiliary approximation result for  conformally invariant functionals at a fixed $u\in\hat{H}^{1}\cap L^{\infty}$ by means of a sequence of functions $u^{n}$ with prescribed compact support. For $G\in C^{0}(\R^{3},\R^{3})$ let us denote
$$
\G(u):=\int_{\R^{2}}G(u)\cdot u_{x}\wedge u_{y}\quad (u\in\hat{H}^{1}\cap L^{\infty}).
$$

\begin{Lemma}
\label{L:auxiliary}
Let $D$ be an open disc in $\R^{2}$. 
\begin{itemize}[leftmargin=20pt]
\item[(i)] 
For every $u\in\hat{H}^{1}\cap L^{\infty}$ with $u(z)=p\in\R^{3}$ for $|z|$ large, there exists a sequence $(u^{n})\subset H^{1}_{0}(D,\R^{3})\cap L^{\infty}$ such that 
\begin{gather}
\label{eq:auxiliary-1}
\|u^{n}\|_{\infty}\le\|u\|_{\infty},\quad\D(u^{n})\to\D(u),\\
\label{eq:auxiliary-2}
\G(u^{n})=\G(u)\quad\forall G\in C^{0}(\R^{3},\R^{3}).
\end{gather}
If in addition $u\in\R^{3}+C^{\infty}_{c}(\R^{2},\R^{3})$, then the sequence $(u^{n})$ can be taken in $C^{\infty}_{c}(D,\R^{3})$. 
\item[(ii)] 
For every $u\in\hat{H}^{1}\cap L^{\infty}$ with $u(z)\to p\in\R^{3}$ as $|z|\to\infty$ and $\V(u)\ne 0$, there exists a sequence $(u^{n})\subset H^{1}_{0}(D,\R^{3})\cap L^{\infty}$ satisfying (\ref{eq:auxiliary-1}), $\V(u^{n})=\V(u)$, and $\G(u^{n})\to\G(u)$ for finitely many $G\in C^{0}(\R^{3},\R^{3})$.
\end{itemize}
\end{Lemma}

\Proof
(i) Let $u\in\hat{H}^{1}\cap L^{\infty}$ with $u(z)=p\in\R^{3}$ for $|z|\ge R$. For every integer $n>\max\{1,R\}$ let $\eta_{n}\colon\R\to[0,1]$ be a smooth decreasing function such that
$$
\eta_{n}(r)=\left\{\begin{array}{ll}1&\text{as }r\le n\\ 2-\frac{\log r}{\log n}&\text{as }n+1<r\le n^{2}-1\\ 0&\text{as }r>n^{2}.\end{array}\right.
$$
Then set
\begin{equation}
\label{eq:tilde-un}
\tilde{u}^{n}(z)=\left\{\begin{array}{ll}u(z)&\text{as }|z|\le n\\ \eta_{n}(|z|)p&\text{as }n<|z|\le n^{2}\\ 0&\text{as }|z|>n^{2}.\end{array}\right.
\end{equation}
By direct computations, one can check that (\ref{eq:auxiliary-1}) and (\ref{eq:auxiliary-2}) hold true for $(\tilde{u}^{n})$. Notice that $\tilde{u}^{n}\in H^{1}_{0}(\Omega_{n},\R^{3})$, where $\Omega_{n}$ denotes the disc centered at the origin and with radius $n^{2}$. Let $D$ be a disc centered at some $z_{0}$ and with radius $r>0$. Setting $u^{n}(z)=\tilde{u}^{n}\left(\frac{n^{2}}{r}(z-z_{0})\right)$, one has that $u^{n}\in H^{1}_{0}(D,\R^{3})$ with $\|u^{n}\|_{\infty}=\|\tilde{u}^{n}\|_{\infty}$, $\D(u^{n})=\D(\tilde{u}^{n})$ and $\G(u^{n})=\G(\tilde{u}^{n})$. Hence (\ref{eq:auxiliary-1}) and (\ref{eq:auxiliary-2}) hold true also for $({u}^{n})$. Moreover, if $u$ is smooth, then according to the definition (\ref{eq:tilde-un}), also $\tilde{u}^{n}$ and consequently $u^{n}$ are so.\\
(ii) Let $u\in\hat{H}^{1}\cap L^{\infty}$ with $u(z)\to p\in\R^{3}$ as $|z|\to\infty$. Consider the sequence $(u^{n})$ defined by (\ref{eq:approx-u}). Then, following the proof of Lemma \ref{L:density}, one can recognize that 
\begin{equation}
\label{eq:step0}
u^{n}\to u\text{ in $\hat{H}^{1}$ and in $L^{\infty}_{loc}$~\!, and $\|u^{n}\|_{\infty}\le\|u\|_{\infty}$.}
\end{equation} 
Moreover $u^{n}(z)=p$ for $|z|\ge n^{2}$. Hence we are in the position to apply part (i): for every $n\in\N$ there is a sequence $(\tilde{u}^{n,k})_{k>n}\subset\hat{H}^{1}$ such that 
\begin{equation}
\label{eq:step1}
\begin{array}{c}
\tilde{u}^{n,k}\in H^{1}_{0}(\Omega_{n},\R^{3}),~~\|\tilde{u}^{n,k}\|_{\infty}\le\|u^{n}\|_{\infty},~~\D(\tilde{u}^{n,k})\to\D(u^{n})\text{ as }k\to\infty,\\
\G(\tilde{u}^{n,k})=\G(u^{n})\text{ for every $G\in C^{0}(\R^{3},\R^{3})$ and for $k>n$.}\end{array}
\end{equation}
Let $(\varepsilon_{h})\subset(0,\infty)$ be a sequence such that $\varepsilon_{h}\to 0$. Then, there exists a sequence $n_{h}\to\infty$ such that for every $h\in\N$ one has 
\begin{equation}
\label{eq:step2}
\begin{array}{c}
|\D(u^{n_{h}})-\D(u)|<\varepsilon_{h},~~|\V(u^{n_{h}})-\V(u)|<\varepsilon_{h},\\
|\G(u^{n_{h}})-\G(u)|<\varepsilon_{h}\text{ for finitely many vector fields $G$ in a fixed set $\mathscr{G}$.}
\end{array}
\end{equation} 
In particular the last inequality is justified as follows: for (\ref{eq:step0}), one has that $u^{n}_{x}\wedge u^{n}_{y}\to u_{x}\wedge u_{y}$ in $L^{1}(\R^{2},\R^{3})$ and $\sup_{n}\|G\circ u^{n}\|_{\infty}<\infty$. Hence the dominated convergence theorem applies and one can infer that $\G(u^{n})\to\G(u)$. From (\ref{eq:step1})--(\ref{eq:step2}), for every $h\in\N$ one can find $k_{h}>n_{h}$ such that $|\D(u^{n_{h}})-\D(\tilde{u}^{n_{h},k_{h}})|<\varepsilon_{h}$. Moreover we have $\|\tilde{u}^{n_{h},k_{h}}\|_{\infty}\le\|u^{n_{h}}\|_{\infty}$ and $\G(\tilde{u}^{n_{h},k_{h}})=\G(u^{n_{h}})$ for every $G\in C^{0}(\R^{3},\R^{3})$. Hence, setting $\tilde{u}^{h}=\tilde{u}^{n_{h},k_{h}}$, one has that 
\begin{gather*}
\tilde{u}^{h}\in H^{1}_{0}(\Omega_{n_{h}},\R^{3}),\quad\|\tilde{u}^{h}\|_{\infty}\le\|u\|_{\infty},\quad\D(\tilde{u}^{h})\to\D(u),\\
\V(\tilde{u}^{h})\to\V(u),\quad\G(\tilde{u}^{h})\to\G(u)\quad\forall G\in\mathscr{G}.
\end{gather*}
Recalling that, as in part (i), $D$ is the disc centered at $z_{0}$ and with radius $r$, and setting $v^{h}(z)=\tilde{u}^{h}\left(\frac{k_{h}^{2}}{r}(z-z_{0})\right)$, one has that 
\begin{gather*}
v^{h}\in H^{1}_{0}(D,\R^{3}),\quad\|v^{h}\|_{\infty}\le\|u\|_{\infty},\quad\D(v^{h})=\D(\tilde{u}^{h}),\\
\V(v^{h})=\V(\tilde{u}^{h}),\quad\G(v^{h})=\G(\tilde{u}^{h})\quad\forall G\in\mathscr{G}.
\end{gather*}
Finally we normalize each $v^{h}$ in order to fix the volume. To this extent, let $s_{h}=\sqrt[3]{\V(u)/\V(v^{h})}$ and $w^{h}=s_{h}v^{h}$. Then $s_{h}\to 1$, $w^{h}\in H^{1}_{0}(D,\R^{3})$, $\V(w^{h})=\V(u)$, $\D(w^{h})=s_{h}^{2}\D(v^{h})\to\D(u)$, and 
$$
\G(w^{h})=\G(s_{h}\tilde{u}^{h})=s_{h}^{2}\int_{\R^{2}}G(s_{h}\tilde{u}^{h})\cdot\tilde{u}^{h}_{x}\wedge\tilde{u}^{h}_{y}\to\int_{\R^{2}}G(u)\cdot u_{x}\wedge u_{y}.
$$
Indeed $\tilde{u}^{h}_{x}\wedge\tilde{u}^{h}_{y}\to u_{x}\wedge u_{y}$ in $L^{1}(\R^{2},\R^{3})$ and $\sup_{h}\|G\circ(s_{h}\tilde{u}^{h})\|_{\infty}<\infty$, because $\|s_{h}\tilde{u}^{h}\|_{\infty}\le(1+o(1))\|u\|_{\infty}$. Moreover  $G\circ(s_{h}\tilde{u}^{h})\to G\circ u$ pointwise a.e., because $s_{h}\to 1$, $\tilde{u}^{h}(z)=u^{n_{h}}(z)=u(z)$ for $|z|<n_{h}$ and $n_{h}\to\infty$. Hence the sequence $(w^{h})$ satisfies the required properties, and the proof of part (ii) is complete. 
\qed

\section{Isovolumetric problems}

In this section we aim to investigate a family of constrained minimization problems, defined as follows. For every $t\in\R$ set 
\begin{equation}
\label{eq:isovolumetric}
S_{K}(t):=\inf_{u\in M_{t}}\E(u)\quad\text{where}\quad\E(u):=\D(u)+\Q(u)
\end{equation}
and $M_{t}$ is defined in (\ref{eq:Mt-def}). Our ultimate goal is to prove Theorems \ref{T:existence1} and \ref{T:existence2}. Hence, unless differently specified, we always assume that $K\in C^{1}(\R^{3})$ satisfy $(K_{1})$ and $(K_{2})$. The additional assumptions $(K_{3})$ or $(K_{4})$ will be recalled when they will be needed. 

Firstly we point out that the mapping $t\mapsto S_{K}(t)$ is well posed from $\R$ into $\R$, in view of (\ref{eq:Q<1}) and (\ref{eq:VK-estimate}), and can be named the \emph{isovolumetric function}. We also set 
$$
\Q_{t}(u):=\int_{\R^{2}}Q_{K}(tu)\cdot u_{x}\wedge u_{y}\quad(u\in\hat{H}^{1})
$$
and we notice that Remark \ref{R:K-volume} holds true also for the functional $\Q_{t}$. Moreover we introduce the \emph{normalized isovolumetric function} $t\mapsto\tilde{S}_{K}(t)$ defined by
\begin{equation}
\label{eq:StildeK-def}
\tilde{S}_{K}(t):=\inf_{u\in M_{1}}\E_{t}(u)\quad\text{where}\quad\E_{t}(u):=\D(u)+\Q_{t}(u).
\end{equation}
\begin{Remark} 
\label{R:K=0}
For $t=0$ the class $M_{t}$ contains the constant functions. Since $0\le(1-\|Q_{K}\|_{\infty})\D(u)\le\E(u)$, one infers that $S_{K}(0)=0$ and minimizers for $S_{K}(0)$ are exactly the constant functions. Instead, for $t=0$ the vector field $p\mapsto Q_{K}(tp)$ is constant and then $\Q_{t}(u)=0$ for every $u\in\hat{H}^{1}$. Hence, by (\ref{eq:isoperimetric}), $\tilde{S}_{K}(0)=\inf\{\D(u)~|~u\in M_{1}\}=S=\sqrt[3]{36\pi}$, the isoperimetric constant.\\
Let us examine the case $K=0$ and $t\in\R$ fixed. Then $\E=\D$ and, by (\ref{eq:isoperimetric}), $S_{0}(t)=\inf\{\D(u)~|~u\in M_{t}\}=St^{2/3}$. Instead $\tilde{S}_{0}(t)=S$ for every $t\in\R$. 
\end{Remark}

Let us state some preliminary properties of the isovolumetric function $S_{K}(t)$. 

\begin{Lemma}
\label{L:SK(t)}
For every $t\in\R$ the following facts hold:
\begin{itemize}[leftmargin=22pt]
\item[(i)] $S_{K}(-t)=S_{-K}(t)$;
\item[(ii)] $S_{K}(t)=t^{2/3}\tilde{S}_{K}(t^{1/3})$;
\item[(iii)] $S_{K}(t)=S_{K(\cdot+p)}(t)$ for every $p\in\R^{3}$.
\item[(iv)] $S_{K}(t)=\inf\{\E(u)~|~u\in C^{\infty}_{c}(\R^{2},\R^{3}),~\V(u)=t\}$.
\end{itemize}
\end{Lemma}

\Proof
(i) For every $u\in\hat{H}^{1}$ let $\overline{u}(x,y)=u(y,x)$. Then $\D(\overline{u})=\D(u)$, $\V(\overline{u})=-\V(u)$ and $\Q(\overline{u})=-\Q(u)=\int_{\R^{2}}Q_{-K}(u)\cdot u_{x}\wedge u_{y}$. These identities easily imply that $S_{K}(-t)=S_{-K}(t)$.\\
(ii) For every $u\in\hat{H}^{1}$ and $t\in\R$ one has that 
$\Q(tu)=t^{2}\Q_{t}(u)$ whereas $\D(tu)=t^{2}\D(u)$ and $\V(tu)=t^{3}\V(u)$ and thus $S_{K}(t)=t^{2/3}\tilde{S}_{K}(t^{1/3})$.\\
(iii) Fix $p\in\R^{3}$. For every $u\in\hat{H}^{1}$ one has that $\D(u+p)=\D(u)$, $\V(u+p)=\V(u)$, and $\Q(u+p)=\int_{\R^{2}}Q_{K(\cdot+p)}(u)\cdot u_{x}\wedge u_{y}$. Then 
$$
S_{K(\cdot+p)}(t)=\inf\{\E(u)~|~u\in p+\hat{H}^{1},~\V(u)=t\}=S_{K}(t)
$$
because $p+\hat{H}^{1}=\hat{H}^{1}$.\\
(iv) Fix $\varepsilon>0$ and take $u\in\hat{H}^{1}$ with $\V(u)=t$  and $\E(u)\le S_{K}(t)+\varepsilon$. By Lemma \ref{L:density} and by the continuity of the functionals $\D$, $\V$, and $\Q$, there exists a sequence $(u^{n})\subset\R^{3}+C^{\infty}_{c}(\R^{2},\R^{3})$ such that 
$$
\E(u^{n})\le S_{K}(t)+\varepsilon+o(1)\quad\text{and}\quad\V(u^{n})=t+o(1),
$$ 
where $o(1)\to 0$ as $n\to\infty$. We can write $u^{n}=p^{n}+\tilde{u}^{n}$ with $p^{n}\in\R^{3}$ and $\tilde{u}^{n}\in C^{\infty}_{c}(\R^{2},\R^{3})$. Let $w\in C^{\infty}_{c}(\R^{2},\R^{3})$ be a mapping with $\V(w)=1$. We can find a sequence $z_{n}\in\R^{2}$ such that $w^{n}=w(\cdot +z_{n})$ has support with empty intersection with the support of $\tilde{u}^{n}$. Notice that $\D(w^{n})=\D(w)$ and $\V(w^{n})=1$ for all $n\in\N$. Finally, we define
$$
v^{n}=u^{n}+s_{n}w^{n}\quad\text{with}\quad s_{n}=\sqrt[3]{t-\V(u^{n})}.
$$
We have that $v^{n}\in p_{n}+C^{\infty}_{c}(\R^{2},\R^{3})$. Moreover, using (\ref{eq:VK-estimate}) and since $s_{n}\to 0$ as $n\to\infty$, we estimate
\begin{equation*}
\begin{split}
&\V(v^{n})=\V(u^{n})+s_{n}^{3}\V(w^{n})=t\\
&\D(v^{n})=\D(u^{n})+s_{n}^{2}\D(w)=\D(u^{n})+o(1)\\
&\Q(v^{n})=\Q(u^{n})+\Q(s_{n}w^{n})=\Q(u^{n})+o(1)
\end{split}
\end{equation*} 
where $o(1)\to 0$ as $n\to\infty$. Hence for fixed $n$ large enough, 
\begin{equation}
\label{eq:pn}
\inf\{\E(u)~|~u\in p_{n}+C^{\infty}_{c}(\R^{2},\R^{3}),~\V(u)=t\}\le\E(v^{n})\le S_{K}(t)+2\varepsilon.
\end{equation}
Now we claim that for every $p\in\R^{3}$
\begin{equation}
\label{eq:claim-p}
\inf_{M_{t}\cap  C^{\infty}_{c}(\R^{2},\R^{3})}\E\le\tilde{S}_{K,p}(t):=\inf\{\E(u)~|~u\in p+C^{\infty}_{c}(\R^{2},\R^{3}),~\V(u)=t\}.
\end{equation}
Indeed, fixing $\varepsilon>0$, let $u\in C^{\infty}(\R^{2},\R^{3})$ be such that $u(z)=p$ for $|z|\ge R$, $\V(u)=t$ and 
$$
\E(u)\le\tilde{S}_{K,p}(t)+\varepsilon.
$$
By Lemma \ref{L:auxiliary} (i), there exists a sequence $(u^{n})\subset C^{\infty}_{c}(\R^{2},\R^{3})$ such that $\D(u^{n})\to\D(u)$, $\V(u^{n})=\V(u)$ and $\Q(u^{n})=\Q(u)$. Hence $\E(u^{n})\to\E(u)$ and for $n$ large enough
$$
\inf_{M_{t}\cap  C^{\infty}_{c}(\R^{2},\R^{3})}\E\le\E(u^{n})\le\tilde{S}_{K,p}(t)+2\varepsilon.
$$ 
Therefore by the arbitrariness of $\varepsilon>0$, (\ref{eq:claim-p}) follows. By (\ref{eq:pn}) and (\ref{eq:claim-p}), we conclude that 
$$
\inf_{M_{t}\cap  C^{\infty}_{c}(\R^{2},\R^{3})}\E\le S_{K}(t).
$$
Since the opposite inequality is trivial, (v) is proved.
\qed\medskip

Now we state some estimates on the isovolumetric function $S_{K}(t)$.

\begin{Lemma}
\label{L:SK-estimates}
For every $t\in\R$ the following facts hold:
\begin{itemize}[leftmargin=20pt]
\item[(i)] For every $t\in\R$ one has that $(1-\|Q_{K}\|_{\infty})St^{2/3}\le S_{K}(t)\le S_{0}(t)=St^{2/3}$.
\item[(ii)] For every $t_{1},...,t_{k}\in\R$ one has that $S_{K}(t_{1})+...+S_{K}(t_{k})\ge S_{K}(t_{1}+...+t_{k})$.
\end{itemize}
\end{Lemma}

\Proof
(i) The first inequality follows from (\ref{eq:isoperimetric}) and (\ref{eq:VK-estimate}). Let us show the second one. Since $K$ satisfies $(K_{1})$--$(K_{2})$ if and only if $-K$ does so, by Lemma \ref{L:SK(t)} (i), without loss of generality we can assume $t<0$. Let $p^{n}\in\R^{3}$ be such that $|p^{n}|\to\infty$ and let $u^{n}=r\omega+p^{n}$ where $\omega$ is defined in Remark \ref{R:sphere} and $r>0$ is such that $-4\pi r^{3}/3=t$. Then $u^{n}\in M_{t}$ and $\E(u^{n})=t^{2/3}S-\int_{B_{r}(p^{n})}K(p)~\!dp$ (see Remarks \ref{R:sphere} and \ref{R:sphere2}) and the conclusion follows from the fact that, by $(K_{2})$, $K(p)\to 0$ as $|p|\to\infty$. 
\\
(ii) Let $t_{1},...,t_{k}\in\R$ be given and fix an arbitrary $\varepsilon>0$. By Lemma \ref{L:SK(t)} (iv) there exist $u^{1},...,u^{k}\in C^{\infty}_{c}(\R^{2},\R^{3})$ such that
$$
\V(u^{i})=t_{i},~~\D(u^{i})+\Q(u^{i})\le S_{K}(t_{i})+\frac{\varepsilon}{k}\quad\forall i=1,...,k.
$$
Up to translation we can assume that the supports of the mappings $u^{i}$ are pairwise disjoint. Then $\V\left(\sum_{i}u^{i}\right)=\sum_{i}t_{i}$ and $S_{K}\textstyle{\left(\sum_{i}t_{i}\right)}\le \E\textstyle{\left(\sum_{i}u^{i}\right)}=\sum_{i}\E(u^{i})\le \sum_{i}S_{K}(t_{i})+\varepsilon$. By the arbitrariness of $\varepsilon>0$, (ii) holds. 
\qed\medskip

The next result contains some properties about minimizing sequences for the isovolumetric problem defined by (\ref{eq:isovolumetric}). In particular we state a bound from above and from below on the Dirichlet norm, and we show that every minimizing sequence shadows another minimizing sequence consisting of approximating solutions for some prescribed mean curvature equation. 

\begin{Lemma}
\label{L:minimizing-PS}
Let $t\in\R$ be fixed. Then:
\begin{itemize}[leftmargin=22pt]
\item[(i)] 
$\D(u)\ge \frac{S_{K}(t)}{1+\|Q_{K}\|_{\infty}}$ for every $u\in M_{t}$.
\item[(ii)]
If $(u^{n})\subset M_{t}$ is a minimizing sequence for $S_{K}(t)$ then $\limsup\D(u^{n})\le \frac{St^{2/3}}{1-\|Q_{K}\|_{\infty}}$.
\item[(iii)]
For every minimizing sequence $(\tilde{u}^{n})\subset M_{t}$ for $S_{K}(t)$ there exists another minimizing sequence $(u^{n})\subset M_{t}$ such that $\|u^{n}-\tilde{u}^{n}\|\to 0$ and with the additional property that
\begin{equation}
\label{eq:almost-solution}
\Delta u^{n}-K(u^{n})u^{n}_{x}\wedge u^{n}_{y}+\lambda u^{n}_{x}\wedge u^{n}_{y}\to 0\quad\text{in }\hat{H}^{-1}\text{(= dual of $\hat{H}^{1}$)}
\end{equation}
for some $\lambda\in\R$. 
\end{itemize}
\end{Lemma}

\Proof
(i) and (ii) can be easily obtained by (\ref{eq:VK-estimate}) and by Lemma \ref{L:SK-estimates} (i).\\
(iii) Assume $t=0$. Then $S_{K}(0)=0$ (see Remark \ref{R:K=0}) and if $(\tilde{u}^{n})\subset M_{0}$ is a minimizing sequence for $S_{K}(0)$ then $\D(\tilde{u}^{n})\to 0$ by part (ii). Taking $u^{n}=\frac{1}{4\pi}\int_{\R^{2}}\tilde{u}^{n}\mu^{2}$, one easily checks that $(u^{n})$ satisfies the thesis by the Poincar\'e inequality which leads to $\|u^{n}-\tilde{u}^{n}\|\to 0$ as $n\to\infty$. (Indeed each $u^{n}$ is a constant and is a minimizer for $S_{K}(0)$). Now let us examine the case $t\ne 0$. Since, in general, the functional $\E$ is not differentiable everywhere in $M_{t}$, the proof of (iii) needs some care. Since $t\ne 0$, the set $M_{t}$ constitutes a smooth closed manifold (see Lemma \ref{L:volume}). Let $(\tilde{u}^{n})\subset M_{t}$ be such that $\E(\tilde{u}^{n})\to S_{K}(t)$ and fix a sequence $(\varepsilon_{n})\subset(0,\infty)$ with $\varepsilon_{n}\to 0$. By Ekeland's variational principle (see, e.g., \cite{Eke79}), there exists a sequence $(u^{n})\subset M_{t}$ such that 
$$
\|u^{n}-\tilde{u}^{n}\|\le\varepsilon_{n},\quad\E(u^{n})\le\E(\tilde{u}^{n}),\quad\E(u^{n})\le\E(u)+\varepsilon_{n}\|u-u^{n}\|~~\forall u\in M_{t}.
$$
Fix $\varphi\in T_{u^{n}}M_{t}\cap L^{\infty}$ and for $s>0$ small enough, set $\tau_{n}(s)=\sqrt[3]{t/\V(u^{n}+s\varphi)}$. Then $\tau_{n}(s)(u^{n}+s\varphi)\in M_{t}$ and 
\begin{equation}
\label{eq:Ekeland-inequality}
\frac{\E(\tau_{n}(s)(u^{n}+s\varphi))-\E(u^{n})}{s}\ge-\varepsilon_{n}\left\|\frac{\tau_{n}(s)(u^{n}+s\varphi)-u^{n}}{s}\right\|.
\end{equation}
We compute the limit as $s\to 0^{+}$ in the following separate auxiliary Lemma.

\begin{Lemma}
\label{L:technical}
Let $t\in\R\setminus\{0\}$ be fixed. For every $u\in M_{t}$ and $\varphi\in\hat{H}^{1}\cap L^{\infty}$ with $\V'(u)[\varphi]=0$, it holds that:
\begin{gather}
\label{eq:technical-1}
\lim_{s\to 0^{+}}\left\|\frac{\tau(s)(u+s\varphi)-u}{s}-\varphi\right\|=0\\
\label{eq:technical-2}
\lim_{s\to 0^{+}}\frac{\Q(\tau(s)(u+s\varphi))-\Q(u)}{s}=\int_{\R^{2}}K(u)\varphi\cdot u_{x}\wedge u_{y}\\
\label{eq:technical-3}
\lim_{s\to 0^{+}}\frac{\E(\tau(s)(u+s\varphi))-\E(u)}{s}=\int_{\R^{2}}(\nabla u\cdot\nabla\varphi+K(u)\varphi\cdot u_{x}\wedge u_{y})=:\E'(u)[\varphi].
\end{gather}
where $\tau(s)=\sqrt[3]{t/\V(u+s\varphi)}$. 
\end{Lemma}
Hence, passing to the limit as $s\to 0^{+}$ in (\ref{eq:Ekeland-inequality}), by Lemma \ref{L:technical} we obtain that $\E'(u^{n})[\varphi]\ge-\varepsilon_{n}\|\varphi\|$. Taking now $-\varphi$ instead of $\varphi$ we get $\E'(u^{n})[\varphi]\le\varepsilon_{n}\|\varphi\|$. Then, since $\hat{H}^{1}\cap L^{\infty}$ is dense in $\hat{H}^{1}$ we conclude that 
\begin{equation}
\label{eq:Ekeland-inequality-2}
\sup_{\scriptstyle\varphi\in T_{u^{n}}M_{t}\atop\scriptstyle\varphi\ne 0}\frac{|\E'(u^{n})[\varphi]|}{\|\varphi\|}\le\varepsilon_{n}.
\end{equation}
Now let $v^{n}\in\hat{H}^{1}$ be the Riesz representative of $\V'(u^{n})$. Set 
$$
\lambda_{n}=\frac{\E'(u^{n})[v^{n}]}{\|v^{n}\|^{2}}
$$
(notice that $\lambda_{n}$ is well defined because $v^{n}\in L^{\infty}$, see Lemma \ref{L:volume}). For every $\varphi\in\hat{H}^{1}\cap L^{\infty}$ the projection of $\varphi$ on $T_{u^{n}}M_{t}$ is given by
$$
\tilde\varphi=\varphi-\frac{\langle v^{n},\varphi\rangle}{\|v^{n}\|^{2}}v_{n}
$$
and, by (\ref{eq:Ekeland-inequality-2}), 
$$
|\E'(u^{n})[\varphi]-\lambda_{n}\V'(u^{n})[\varphi]|=|\E'(u^{n})[\tilde\varphi]|\le\varepsilon_{n}\|\tilde\varphi\|\le\varepsilon_{n}\|\varphi\|,
$$
and then, by density, $\E'(u^{n})-\lambda_{n}\V'(u^{n})\to 0$ in $\hat{H}^{-1}$. Now we show that the sequence $(\lambda_{n})$ is bounded. Indeed, by (\ref{eq:Wente-inequality}) and by Lemma  \ref{L:minimizing-PS}, part (ii), we know that 
\begin{equation}
\label{eq:vn-bounded}
\|\nabla v^{n}\|_{2}+\|v^{n}\|_{\infty}\le C\|\nabla u^{n}\|_{2}^{2}\le C.
\end{equation} 
Then
\begin{eqnarray}
\nonumber
|\E'(u^{n})[v^{n}]|&\le&\left|\int_{\R^{2}}(\nabla u^{n}\cdot\nabla v^{n}+K(u^{n})v^{n}\cdot u^{n}_{x}\wedge u^{n}_{y})\right|\\
\label{eq:lambdan-bdd-1}
&\le&\|\nabla u^{n}\|_{2}\|\nabla v^{n}\|_{2}+\|K\|_{\infty}\|v^{n}\|_{\infty}\|\nabla u^{n}\|_{2}^{2}\le C.
\end{eqnarray}
Moreover, keeping into account that $\int_{\R^{2}}v^{n}\mu^{2}=0$ and using again Lemma  \ref{L:minimizing-PS}, part (ii), we have that
\begin{eqnarray}
\nonumber
|3t|&=&|\V'(u^{n})[u^{n}]|=|\langle v^{n},u^{n}\rangle|=\left|\int_{\R^{2}}\nabla v^{n}\cdot\nabla u^{n}\right|\\
\label{eq:lambdan-bdd-2}
&\le&\|\nabla v^{n}\|_{2}\|\nabla u^{n}\|_{2}\le C\|\nabla v^{n}\|_{2}=C\|v^{n}\|.
\end{eqnarray}
Then (\ref{eq:lambdan-bdd-1}) and (\ref{eq:lambdan-bdd-2}) imply that $(\lambda_{n})$ is bounded, because $t\ne 0$. Hence, for a subsequence $\lambda_{n}\to\lambda\in\R$ and since $(v^{n})$ is bounded in $\hat{H}^{1}$ (use (\ref{eq:vn-bounded})), we conclude that $\E(u^{n})-\lambda\V'(u^{n})\to 0$ in $\hat{H}^{-1}.\quad\square$
\medskip

\noindent
\emph{Proof of Lemma \ref{L:technical}.}
First of all we observe that, by Lemma \ref{L:volume}, the mapping $s\mapsto\tau(s)$ is smooth, with 
$$
\tau'(s)=-\frac{\sqrt[3]{t}}{3}\frac{\V'(u+s\varphi)[\varphi]}{\V(u+s\varphi)^{4/3}}.
$$
In particular, 
\begin{equation}
\label{eq:technical-4}
\tau(0)=1,\quad\lim_{s\to 0}\frac{\tau(s)-1}{s}=\tau'(0)=0
\end{equation}
because $\V'(u)[\varphi]=0$. Hence (\ref{eq:technical-1}) easily follows from (\ref{eq:technical-4}). In order to prove (\ref{eq:technical-2}) we write
$$
\frac{\Q(\tau(s)(u+s\varphi))-\Q(u)}{s}=
\frac{\Q(\tau(s)(u+s\varphi))-\Q(\tau(s)u)}{s}+\frac{\Q(\tau(s)u)-\Q(u)}{s}.
$$
Using (\ref{eq:radial-VK-derivative}), we have that 
$$
\lim_{\tau\to 1}\frac{\Q(\tau u)-\Q(u)}{\tau-1}=\int_{\R^{2}}K(u)u\cdot u_{x}\wedge u_{y}\in\R
$$
and consequently
\begin{equation}
\label{eq:technical-5}
\lim_{s\to 0}\frac{\Q(\tau(s)u)-\Q(u)}{s}=\lim_{s\to 0}\frac{\tau(s)-1}{s}~\!\lim_{\tau\to 1}\frac{\Q(\tau u)-\Q(u)}{\tau-1}=0.
\end{equation}
Then, writing $u^{s}=\tau(s)u$ and $\varphi^{s}=\tau(s)\varphi$, by Lemma \ref{L:K-volume}, part (ii),
\begin{equation*}\begin{split}
&\frac{\Q(\tau(s)(u+s\varphi))-\Q(\tau(s)u)}{s}\\
&\qquad=\int_{\R^{2}}\int_{0}^{1}K(u^{s}+rs\varphi^{s})\varphi^{s}\cdot(u^{s}_{x}+rs\varphi^{s}_{x})\wedge(u^{s}_{y}+rs\varphi^{s}_{y})~\!dr~\!dz.
\end{split}\end{equation*}
We point out that if $s\to 0$ then $u^{s}\to u$ and $\varphi^{s}\to\varphi$ in $\hat{H}^{1}$ and pointwise a.e. Hence 
\begin{gather*}
(u^{s}_{x}+rs\varphi^{s}_{x})\wedge(u^{s}_{y}+rs\varphi^{s}_{y})\to u_{x}\wedge u_{y}\quad\text{in }L^{1},\\
K(u^{s}+rs\varphi^{s})\to K(u)\quad\text{pointwise a.e.} 
\end{gather*}
and $\|K(u^{s}+rs\varphi^{s})\varphi^{s}\|_{\infty}\le C$ for $s$ close to $0$, because $K$ and $\varphi$ are bounded functions. Then
\begin{equation}
\label{eq:technical-6}
\lim_{s\to 0}\frac{\Q(\tau(s)(u+s\varphi))-\Q(\tau(s)u)}{s}=\int_{\R^{2}}K(u)\varphi\cdot u_{x}\wedge u_{y}
\end{equation}
and (\ref{eq:technical-2}) follows from (\ref{eq:technical-5}) and (\ref{eq:technical-6}). Finally, since $\E=\D+\Q$,  (\ref{eq:technical-3}) is an obvious consequence of (\ref{eq:technical-2}) and of the fact that the functional $\D$ is analytic. 
\qed\medskip

As a next step, we provide a precise description of the specific minimizing sequences for the isovolumetric problems (\ref{eq:isovolumetric}). To this purpose, it is convenient to introduce the definition of a bubble: 

\begin{Definition}
Let $H\in C^{0}(\R^{3})$ be a given function. We call $U\in\hat{H}^{1}$ an $H$-bubble if it is a nonconstant solution to
\begin{equation}
\label{eq:Hsystem}
\Delta U=H(U)U_{x}\wedge U_{y}\quad\text{on}\quad\R^{2}
\end{equation} 
in the distributional sense. If $H$ is constant, an $H$-bubble will be named $H$-sphere. The system (\ref{eq:Hsystem}) is called $H$-system. 
\end{Definition}

Let us state a first preliminary property of $H$-bubbles, for a class of mappings $H$ of our interest.

\begin{Lemma}
\label{L:bubble-bdd}
Let $H(p)=K(p)-\lambda$ with $\lambda\in\R$ and $K\in C^{0}(\R^{3})$ satisfying $(K_{1})$. If $U\in\hat{H}^{1}$ is an $H$-bubble, then $U\in L^{\infty}$, and $\lambda\V(U)>0$. If, in addition, $K\in C^1(\R^3)$ then 
$U$ is of class $C^{2,\alpha}$ as a map on $\S^{2}$.  
\end{Lemma}

\Proof
Multiplying $\Delta U=H(U)U_{x}\wedge U_{y}$ by $U$ and integrating, by $(K_{1})$, one obtains that
$$
0<(2-k_{0})\D(U)\le 3\lambda\V(U)
$$
because $U$ is nonconstant. The fact that $U\in L^{\infty}$ has been proved in \cite{CaMuJFA}, Theorem 2.1 and Remark 2.5. When $K$, and then $H$, is of class $C^1$, the regularity theory for $H$-systems (\ref{eq:Hsystem}) applies (see, for instance, \cite{Hei86} or \cite{BetGhi93}) and one infers that $U\in C^{2,\alpha}(\R^{2},\R^{3})$. By the invariance of $H$-systems and of the Dirichlet integral with respect to the transformation $(x,y)\mapsto(\hat{x},\hat{y}):=(\frac{x}{x^2+y^2},-\frac{y}{x^2+y^2})$, also the mapping $\hat{U}(x,y):=U(\hat{x},\hat{y})$ is an $H$-bubble. From this one infers that $U$ is of class $C^{2,\alpha}$ as a map on $\S^{2}$.
\qed\medskip

According to the following crucial result, minimizing sequences for problems (\ref{eq:isovolumetric}) admit a limit configuration made by bubbles. More precisely:

\begin{Lemma}[Decomposition Theorem]
\label{L:PS-decomposition}
Let $K\colon\R^{3}\to\R$ be a continuous function satisfying $(K_{1})$ and $(K_{2})$. If $(u^{n})\subset\hat{H}^{1}$ is a sequence  satisfying (\ref{eq:almost-solution}) for some $\lambda\in\R$ and such that $c_{1}\le\|\nabla u^{n}\|_{2}\le c_{2}$ for some $0<c_{1}\le c_{2}<\infty$ and for every $n$, then there exist a subsequence of $(u^{n})$, still denoted $(u^{n})$, finitely many $(K-\lambda)$-bubbles $U^{i}$ ($i\in I$), finitely many $(-\lambda)$-spheres $U^{j}$ ($j\in J$) 
such that, as $n\to\infty$:
\begin{gather}
\label{eq:PS}
\left\{\begin{array}{l}\D(u^{n})\to\textstyle\sum_{i\in I}\D(U^{i})+\sum_{j\in J}\D(U^{j})\\
\V(u^{n})\to\textstyle\sum_{i\in I}\V(U^{i})+\sum_{j\in J}\V(U^{j})\\
\Q(u^{n})\to\textstyle\sum_{i\in I}\Q(U^{i})
\end{array}\right.
\end{gather}
where $I$ or $J$ can be empty but not both. In particular, if $J=\varnothing$ then the subsequence $(u^{n})$ is bounded in $\hat{H}^{1}$. 
\end{Lemma}

\Proof 
This result is obtained by combining Theorem 0.1 of \cite{CaJFA} with the proof of Theorem A.1 of \cite{CaMu11}, which in fact holds true assuming just $(K_{1})$ and $(K_{2})$. See also \cite{CaMuJFA} for previous partial, fundamental results.  
\qed\medskip

The next result will be used to show that minimizing sequences for the isovolumetric problems (\ref{eq:isovolumetric}) do not split in two or more $(K-\lambda)$-bubbles. We point out that assumption $(K_{3})$ plays a role just at this point of the argument. We also recall that, since $(K_{3})$, together with $(K_{2})$, implies $(K_{1})$, Lemma \ref{L:bubble-bdd} applies and, in particular, if $U\in\hat{H}^{1}$ is a $(K-\lambda)$-bubble, then $\lambda\ne0$.  

\begin{Lemma}
\label{L:strict-inequality}
Assume $(K_{2})$--$(K_{3})$. If $U^{1},U^{2}\in\hat{H}^{1}$ are two $(K-\lambda)$-bubbles, for a common $\lambda\in\R$, with $\V(U^{i})=t_{i}>0$ ($i=1,2$), then $S_{K}(t_{1}+t_{2})<\E(U^{1})+\E(U^{2})$.
\end{Lemma}

\Proof
Let $D_{1}$ and $D_{2}$ be two disjoint discs, and for $i=1,2$ let $u^{i}\in H^{1}_{0}(D_{i},\R^{3})$ with $\V(u^{i})=t_{i}$. Then set $\tau=\frac{t_{1}}{t_{2}}$ and for $s\in[0,1+\tau^{-1}]$
$$
v^{s}=\sqrt[3]{s}~\!u^{1}+\sqrt[3]{(1-s)\tau+1}~\!u^{2}.
$$
Note that $\V(v^{s})=s\V(u^{1})+((1-s)\tau+1)\V(u^{2})=t_{1}+t_{2}$. Moreover the mapping 
$$
f(s):=\E(v^{s})
$$
is continuous in $[0,1+\tau^{-1}]$ and since $f(s)=\E(\sqrt[3]{s}~\!u^{1})+\E(\sqrt[3]{(1-s)\tau+1}~\!u^{2})$, by means of Lemma \ref{L:K-volume} (iii), one can compute the derivatives
\begin{equation}
\label{eq:f-prime}
\begin{split}
f'(s)&=\frac{s^{-\frac{1}{3}}}{3}\int_{\R^{2}}(|\nabla u^{1}|^{2}+G_{0}(\sqrt[3]{s}~\!u^{1})\cdot u^{1}_{x}\wedge u^{1}_{y})
\\
&\quad-\frac{\tau[(1-s)\tau+1]^{-\frac{1}{3}}}{3}\int_{\R^{2}}(|\nabla u^{2}|^{2}+G_{0}(\sqrt[3]{(1-s)\tau+1}~\!u^{2})\cdot u^{2}_{x}\wedge u^{2}_{y}),
\end{split}
\end{equation}
\begin{equation*}
\begin{split}
f''(s)&=-\frac{s^{-\frac{4}{3}}}{9}\int_{\R^{2}}(|\nabla u^{1}|^{2}-G_{1}(\sqrt[3]{s}~\!u^{1})\cdot u^{1}_{x}\wedge u^{1}_{y})\\
&\quad-\frac{\tau^{2}((1-s)\tau+1)^{-\frac{4}{3}}}{9}\int_{\R^{2}}(|\nabla u^{2}|^{2}-G_{1}(\sqrt[3]{(1-s)\tau+1}~\!u^{2})\cdot u^{2}_{x}\wedge u^{2}_{y})
\end{split}
\end{equation*}
where $G_{0}(u)=K(u)u$ and $G_{1}(u)=(\nabla K(u)\cdot u)u$. From $(K_{1})$ it follows that
\begin{gather*}
\int_{\R^{2}}|\nabla u^{1}|^{2}+\int_{\R^{2}}G_{0}(\sqrt[3]{s}~\!u^{1})\cdot u^{1}_{x}\wedge u^{1}_{y}\ge (2-k_{0})\D(u^{1})\\
\int_{\R^{2}}|\nabla u^{2}|^{2}+\int_{\R^{2}}G_{0}(\sqrt[3]{(1-s)\tau+1}~\!u^{2})\cdot u^{2}_{x}\wedge u^{2}_{y}\ge (2-k_{0})\D(u^{2}).
\end{gather*}
In particular $f'(s)\to\infty$ as $s\to 0$ whereas $f'(s)\to-\infty$ as $s\to 1+\tau^{-1}$. From $(K_{3})$ it follows that
\begin{gather*}
\int_{\R^{2}}|\nabla u^{1}|^{2}-\int_{\R^{2}}G_{1}(\sqrt[3]{s}~\!u^{1})\cdot u^{1}_{x}\wedge u^{1}_{y}\ge (2-k_{0})\D(u^{1})\\
\int_{\R^{2}}|\nabla u^{2}|^{2}-\int_{\R^{2}}G_{1}(\sqrt[3]{(1-s)\tau+1}~\!u^{2})\cdot u^{2}_{x}\wedge u^{2}_{y}\ge (2-k_{0})\D(u^{2}).
\end{gather*}
In particular $f$ is strictly concave in $[0,1+\tau^{-1}]$ and, by Lemma   \ref{L:SK-estimates} (i) and Lemma \ref{L:minimizing-PS}, there exists a constant $c>0$ independent of $u^{1}$ and $u^{2}$, such that $f''(s)\le -c$ for every $s\in(0,1+\tau^{-1})$. As a consequence, with elementary arguments, one obtains that 
$$
\min\{f(0),f(1+\tau^{-1})\}\le\max_{s\in[0,1+\tau^{-1}]}f(s)-\delta
$$
where $\delta=\frac{c}{8}(1+\tau^{-1})^{2}>0$. Hence, 
\begin{equation}
\label{eq:delta-estimate}
S_{K}(t_{1}+t_{2})\le\max_{s\in[0,1+\tau^{-1}]}\E(\sqrt[3]{s}~\!u^{1}+\sqrt[3]{(1-s)\tau+1}~\!u^{2})-\delta.
\end{equation}
Now let us prove the strict inequality $S_{K}(t_{1}+t_{2})<S_{K}(t_{1})+S_{K}(t_{2})$. For $i=1,2$ let $U^{i}\in M_{t_{i}}$ be $(K-\lambda)$-bubbles. According to Lemma \ref{L:bubble-bdd}, $U^{1}$ and $U^{2}$ are bounded and, since we are assuming $K$ of class $C^1$, again by Lemma \ref{L:bubble-bdd}, there exists $\lim_{|z|\to\infty}U^{i}(z)\in\R^3$ ($i=1,2$). Fixing two disjoint discs $D_{1}$ and $D_{2}$, by Lemma \ref{L:auxiliary} (ii), there exist sequences $(u^{i,n})_{n}\subset H^{1}(D_{i},\R^{3})\cap L^{\infty}$ ($i=1,2$) with $\V(u^{i,n})=t_{i}$, $\D(u^{i,n})\to\D(U^{i})$, $\E(u^{i,n})\to\E(U^{i})$, and 
$$
\int_{\R^{2}}G_{0}(u^{i,n})\cdot u^{i,n}_{x}\wedge u^{i,n}_{y}\to\int_{\R^{2}}G_{0}(U^{i})\cdot U^{i}_{x}\wedge U^{i}_{y}.
$$
For every $n\in\N$ set 
$$
f_{n}(s)=\E(\sqrt[3]{s}~\!u^{1,n}+\sqrt[3]{(1-s)\tau+1}~\!u^{2,n})\quad (s\in[0,1+\tau^{-1}]).
$$
By (\ref{eq:delta-estimate}), for every $n$ we have that 
$$
S_{K}(t_{1}+t_{2})\le f_{n}(s_{n})-\delta
$$
where $s_{n}$ is the (unique) value in $[0,1+\tau^{-1}]$ such that  
$f_{n}(s_{n})=\!\displaystyle\max_{s\in[0,1+\tau^{-1}]}\!f_{n}(s)$. We also compute (see (\ref{eq:f-prime}))
\begin{equation*}
\begin{split}
f'_{n}(1)&=\frac{1}{3}\int_{\R^{2}}\left(|\nabla u^{1,n}|^{2}+G_{0}(u^{1,n})\cdot u^{1,n}_{x}\wedge u^{1,n}_{y}\right)\\
&\quad-\frac{\tau}{3}\int_{\R^{2}}\left(|\nabla u^{2,n}|^{2}+G_{0}(u^{2,n})\cdot u^{2,n}_{x}\wedge u^{2,n}_{y}\right)\\
&=\frac{1}{3}\int_{\R^{2}}\left(|\nabla U^{1}|^{2}+G_{0}(U^{1})\cdot U^{1}_{x}\wedge U^{1}_{y}\right)\\
&\quad-\frac{\tau}{3}\int_{\R^{2}}\left(|\nabla U^{2}|^{2}+G_{0}(U^{2})\cdot U^{2}_{x}\wedge U^{2}_{y}\right)+o(1)\\
&=\lambda\V(U^{1})-\tau\lambda\V(U^{2})+o(1)=o(1)\quad\text{as }n\to\infty
\end{split}
\end{equation*}
where in the last line we exploit the fact that $U^{1}$ and $U^{2}$ are $(K-\lambda)$-bubbles and then
$$
\int_{\R^{2}}(|\nabla U^{i}|^{2}+K(U^{i})U^{i}\cdot U^{i}_{x}\wedge U^{i}_{y})=\lambda\int_{\R^{2}}U^{i}\cdot U^{i}_{x}\wedge U^{i}_{y}=3\lambda\V(U^{i})\quad(i=1,2).
$$
From this and using the fact that $f''_{n}(s)\le-c$ for every $s\in(0,1+\tau^{-1})$ with $c>0$ independent of $n$, we infer that $s_{n}\to 1$. Hence, setting $\rho=\frac{1}{2}\min\{\tau^{-1},1\}$, for $n$ large enough $s_{n}\in [1-\rho,1+\rho]\subset(0,1+\tau^{-1})$ and
$$
|f_{n}(s_{n})-f_{n}(1)|\le|s_{n}-1|\max_{|s-1|\le\rho}|f'_{n}(s)|=o(1)\quad\text{as $n\to\infty$.}
$$
Indeed $\max_{|s-1|\le\rho}|f'_{n}(s)|\le c_{1}$ with $c_{1}$ independent of $n$. To check this, just use (\ref{eq:f-prime}) with $u^{i,n}$ instead of $u^{i}$, and the fact that the sequences $(u^{i,n})_{n}$ are bounded in $\hat{H}^{1}$ and in $L^{\infty}$. Hence
\begin{equation*}\begin{split}
S_{K}(t_{1}+t_{2})&\le f_{n}(1)+o(1)-\delta\\
&=\E(u^{1,n})+\E(u^{2,n})+o(1)-\delta=\E(U^{1})+\E(U^{2})-\delta\quad\text{as $n\to\infty$}
\end{split}
\end{equation*}
and this completes the proof. 
\qed\medskip

Finally we can show the existence of minimizers for problems  (\ref{eq:isovolumetric}). 

\begin{Lemma}
\label{L:minimizer}
Assume $(K_{2})$--$(K_{3})$ or, as an alternative, $(K_{1})$--$(K_{2})$ and $(K_{4})$. If $t>0$ is such that
$$
S_{K}(\tau)<S_{0}(\tau)\quad\forall \tau\in(0,t],
\leqno{(*)}
$$
then there exists $U\in M_{t}$ such that $\E(U)=S_{K}(t)$. Moreover such $U$ is a $(K-\lambda)$-bubble, for some $\lambda>0$.
\end{Lemma}

\begin{Remark}\label{R:star}
Recall that $S_{0}(\tau)$ is the infimum value for the Dirichlet integral in the class $M_{\tau}$ of mappings in $\hat{H}^{1}$ parametrizing surfaces with volume $\tau$. We know that $S_{0}(\tau)$ is attained by a conformal parametrization of a round sphere of volume $\tau$ with arbitrary center (Lemma \ref{L:isoperimetric}). On the other hand, $S_{K}(\tau)$ is is the infimum value for the functional $\E=\D+\Q$ in the same class $M_{\tau}$, and $\Q$ has the meaning of $K$-weighted algebraic volume (see Remark \ref{R:sphere2}; see also \cite{CaMuARMA}, Sect.~2.3). Hence, roughly speaking, the inequality $S_{K}(\tau)<S_{0}(\tau)$ means that one can find a ball $B$ of volume $\tau$, with center possibly depending on $\tau$ and $K$, such that $\int_{B}K(p)~\!dp<0$ (see Lemma \ref{L:strict-inequality-SK}, later on).
\end{Remark}
\noindent
\Proof
By Lemma \ref{L:minimizing-PS} there exists a minimizing sequence $(u^{n})\subset M_{t}$ satisfying (\ref{eq:almost-solution}) for some $\lambda\in\R$. By Lemma \ref{L:PS-decomposition} there exist finitely many $(K-\lambda)$-bubbles $U^{i}\in\hat{H}^{1}$ ($i\in I$) and finitely many $(-\lambda)$-spheres $U^{j}\in\hat{H}^{1}$ ($j\in J$) such that, for a subsequence, (\ref{eq:PS}) holds. Recall that $I$ or $J$ (but not both) can be empty. Thus, setting $t_{i}=\V(U^{i})$ for $i\in I\cup J$,
\begin{equation}
\label{eq:SK-decomposition}
\begin{split}
&t=\textstyle\sum_{i\in I}t_{i}+\textstyle\sum_{j\in J}t_{j}\\ 
&S_{K}(t)=\E(u^{n})+o(1)=\textstyle\sum_{i\in I}\E(U^{i})+\sum_{j\in J}\D(U^{j})+o(1)\quad\text{as }n\to\infty.
\end{split}
\end{equation}
By Lemma \ref{L:bubble-bdd} one has that $t_{i}\lambda>0$ for all $i\in I$ and $t_{j}\lambda>0$ for all $j\in J$. Hence $t_{i},t_{j}\in(0,t]$ for all $i\in I$ and $j\in J$. If $J\ne\varnothing$ then, by $(*)$, (\ref{eq:SK-decomposition}), and by Lemma \ref{L:SK-estimates}, we obtain
\begin{equation*}
\begin{split}
S_{K}(t)&\ge\textstyle\sum_{i\in I}S_{K}(t_{i})+\textstyle\sum_{j\in J}S_{0}(t_{j})\\
&>\textstyle\sum_{i\in I}S_{K}(t_{i})+\sum_{j\in J}S_{K}(t_{j})\ge S_{K}\left(\textstyle\sum_{i\in I}t_{i}+\textstyle\sum_{j\in J}t_{j}\right)=S_{K}(t),
\end{split}
\end{equation*}
a contradiction. Therefore $J=\varnothing$ and, by (\ref{eq:SK-decomposition}) and Lemma \ref{L:SK-estimates} (ii),
\begin{equation}
\label{eq:SK-decomposition2}
0\le\textstyle\sum_{i\in I}(\E(U^{i})-S_{K}(t_{i}))=S_{K}(t)-\textstyle\sum_{i\in I}S_{K}(t_{i})\le 0
\end{equation}
that implies $\E(U^{i})=S_{K}(t_{i})$ for all $i\in I$. Now we claim that $I$ is a singleton. We prove this in two cases, as follows. 
\medskip

%\noindent
\emph{Case 1:} $K$ satisfies $(K_{3})$.
\smallskip

\noindent
If $I$ is not a singleton, by Lemma \ref{L:strict-inequality} and by (\ref{eq:SK-decomposition2}), we reach a contradiction. 
\medskip

%\noindent
\emph{Case 2:} $K$ satisfies $(K_{4})$.
\smallskip

\noindent
Since $J=\varnothing$, by Lemma \ref{L:PS-decomposition}, the sequence $(u^{n})$ is bounded in $\hat{H}^{1}$. Then, testing (\ref{eq:almost-solution}) with $u^{n}$ we have
$$
\int|\nabla u^{n}|^{2}+K(u^{n})u^{n}\cdot u^{n}_{x}\wedge u^{n}_{y}=3\lambda t+o(1)\quad\text{as }n\to\infty
$$
and consequently
\begin{equation*}
%\label{eq:estimate1}
3\lambda t\le(2+k_{0})\D(u^{n})+o(1).
\end{equation*}
Using Lemma \ref{L:minimizing-PS} (ii) and Remark \ref{R:K-volume}, we infer that
\begin{equation}
\label{eq:estimateA}
\frac{S}{3\lambda}\ge\frac{2-k_{0}}{2(2+k_{0})}\sqrt[3]{t}.
\end{equation}
Now assume that there exist at least two $(K-\lambda)$-bubbles $U^{1}$ and $U^{2}$ in the decomposition of $(u^{n})$. From
$$
\int_{\R^{2}}|\nabla U^{i}|^{2}+\int_{\R^{2}}K(U^{i})U^{i}\cdot U^{i}_{x}\wedge U^{i}_{y}=3\lambda t_{i}
$$
it follows that
\begin{equation}
\label{eq:estimateB}
\frac{3\lambda t_{i}}{2-k_{0}}\ge\D(U^{i})\ge St_{i}^{2/3},
\end{equation}
having used $(K_{1})$ and the isoperimetric inequality (\ref{eq:isoperimetric}). Thus (\ref{eq:estimateB}) yields 
\begin{equation}
\label{eq:estimateC}
\frac{S}{3\lambda}\le\frac{\sqrt[3]{t_{i}}}{2-k_{0}}.
\end{equation}
Since $0<t_{1}+t_{2}\le t$, using (\ref{eq:estimateA}) and (\ref{eq:estimateC}) and applying the elemantary estimate $\sqrt[3]{r}+\sqrt[3]{1-r}\le 2^{2/3}$ (take $r=t_{1}/(t_{1}+t_{2})$), we obtain  $2^{2/3}(2+k_{0})\ge(2-k_{0})^{2}$, contrary to $(K_{4})$.
\medskip

\emph{Conclusion.} In both cases, there exists just one $U\in\hat{H}^{1}$ such that $\E(u^{n})\to\E(U)$ and $\V(U)=\V(u^{n})+o(1)=t$. This means that $U$ is a minimizer for $\E$ in $M_{t}$. Moreover $U$ is a $(K-\lambda)$-bubble and, by Lemma \ref{L:bubble-bdd}, $\lambda>0$ because $t>0$.
\qed\medskip

As suggested by Remark \ref{R:star}, condition $(*)$ is connected to the sign of $K$. More precisely:

\begin{Lemma}
\label{L:strict-inequality-SK}
For $\overline{p}\in\R^{3}$ and $r>0$ let $B_{r}(\overline{p})=\{p\in\R^{3}~|~|p-\overline{p}|<r\}$. If $K\le 0$ and $K\not\equiv 0$ in $B_{r}(\overline{p})$, then $S_{K}(t)<S_{0}(t)$ for every $t\in(0,4\pi r^{3}/3]$. In particular, if $K(p)\le 0$ for every $p\in\R^{3}$ and $K\not\equiv 0$, then $S_{K}(t)<S_{0}(t)$ for every $t>0$. 
\end{Lemma}

\begin{Remark}
The global negativeness of $K$ is not a necessary condition to ensure the strict inequality $S_{K}(t)<S_{0}(t)$ for every $t>0$. For example, it is enough that $K$ is negative on the tail of some cone, that is, $K(p)<0$ for every $p=r\sigma$ with $r$ large enough and $\sigma\in\Sigma$ where $\Sigma$ is an open domain in $\S^{2}$.
\end{Remark}

\Proof
Fix $t\in(0,4\pi r^{3}/3]$ and let $\delta>0$ be such that $4\pi\delta^{3}/3=t$. Let $u=\overline{p}-\delta\omega$ where $\omega$ is the standard conformal parametrization of the unit sphere, defined in Remark \ref{R:sphere}. Then $\V(u)=t$ and $\E(u)=4\pi\delta^{2}+\int_{B_{\delta}(\overline{p})}K(p)~\!dp<St^{2/3}=S_{0}(t)$ (see Remarks \ref{R:sphere}, \ref{R:sphere2} and \ref{R:K=0}).
\qed
\medskip

\noindent
\emph{Proof of Theorems \ref{T:existence1} and \ref{T:existence2}.}
Assume that $t_{+}>0$ and set $r_{+}=\sqrt[3]{3t_{+}/4\pi}$. Fix $t\in(0,t_{+})$. Then there exists $\overline{p}\in\R^{3}$, possibly depending on $t$, such that $K\le 0$ and $K\not\equiv 0$ in $B_{r}(\overline{p})$, where $r=\sqrt[3]{3t/4\pi}$. Then, by Lemma \ref{L:strict-inequality-SK}, $S_{K}(\tau)<S_{0}(\tau)$ for every $\tau\in(0,t]$. Therefore we can apply Lemma \ref{L:minimizer} in order to infer that there exists a minimizer $U\in\hat{H}^{1}$ for the minimization problem defined by (\ref{eq:isovolumetric}). Moreover such $U$ is a $(K-\lambda)$-bubble, for some $\lambda>0$. By Lemma \ref{L:bubble-bdd}, $U$ is of class $C^{2,\alpha}$ as a map on $\S^{2}$. Then, with a standard procedure (e.g., considering the weak formulation of (\ref{eq:Hsystem}) and taking variations of the form $U\circ\Phi_{t}$, where $\Phi_{t}$ is a smooth flow on $\S^{2}$), one also infers that $U$ satisfies the conformality conditions. Hence $\A(U)=\D(U)$ and $U$ turns out to be a minimizer also for $\F$, namely is a solution of the original isovolumetric problem (\ref{eq:SKt}). Thus the proof for $t\in(0,t_{+})$ is complete. The case $t=0$ is trivial and already discussed in Remark \ref{R:K=0}. Lastly, if $t_{-}<0$, one can conclude for $t\in(t_{-},0)$ by changing sign to $K$ and using Lemma \ref{L:SK(t)} (i) and what we just proved for $t>0$. 
\qed

\begin{Remark}
The $(K-\lambda)$-bubble $U$ found as minimizer of the isovolumetric problem (\ref{eq:SKt}) describes a parametric surface $S=U(\R^{2}\cup\{\infty\})$ such that $K(p)-\lambda$ equals the mean curvature of $S$ at any regular point $p\in S$ (see, for instance, \cite{CaMuCCM}). In addition, $S$ has at most finitely many branch points (see \cite{GulOssRoy73}). We also notice that $U$ is simple in the sense that it cannot be expressed in the form $U(z)=u(z^{n})$ for some $u\in\hat{H}^{1}$ and $n>1$ integer (here we use complex notation). Indeed, otherwise we should have $\V(u)=t/n$ and $S_{K}(t/n)\le\E(u)=n^{-1}\E(U)$. But in the first part of the proof of Lemma \ref{L:minimizer} we have shown that no decomposition of the form $S_{K}(t)\ge S_{K}(t_{1})+...+S_{K}(t_{n})$ with $t=t_{1}+...+t_{n}$ and $0<t_{i}<t$ ($i=1,...,n$) can occur.
\end{Remark}

Useful bounds on the Lagrange multiplier $\lambda$ can be easily deduced. Indeed, multiplying the system $\Delta U=(K(U)-\lambda)U_{x}\wedge U_{y}$ by $U$, integrating, and exploiting $(K_{1})$, one infers that 
$$
(2-k_{0})\D(U)\le 3\lambda t\le(2+k_{0})\D(U).
$$
Then, by Lemma \ref{L:minimizing-PS}, 
$$
\frac{2(2-k_{0})}{3(2+k_{0})}S_{K}(t)\le\lambda t\le\frac{2(2+k_{0})}{3(2-k_{0})}S_{K}(t)
$$
and finally, by Lemma \ref{L:SK-estimates} (i), for $t>0$ one obtains
\begin{equation}
\label{eq:lambda-bound}
\frac{(2-k_{0})^{2}S}{3(2+k_{0})\sqrt[3]{t}}\le\lambda\le\frac{2(2+k_{0})S}{3(2-k_{0})\sqrt[3]{t}}.
\end{equation}

Let us point out that, as a by-product of Theorem \ref{T:existence2} and with the estimates (\ref{eq:lambda-bound}) we obtain a new existence result for the $H$-bubble problem. This result has a perturbative character, in the same direction of other works like \cite{CaMM}, \cite{CaMuDuke}, \cite{Fel05}, \cite{MuJAM}. 

\begin{Theorem}
\label{T:existence4}
Let $K\in C^{1}(\R^{3})$ satisfy $(K_{1})$, $(K_{2})$, and $(K_{4})$. Then there exists a sequence $(\lambda_{n})\subset\R$ with $|\lambda_{n}|\to\infty$ such that for every $n$ there exists a $(K-\lambda_{n})$-bubble.
\end{Theorem}

\Proof
If $K\equiv 0$ then the result is trivial. If $K\not\equiv 0$ then $(t_{-},t_{+})$ is nonempty, with $t_{-}\le 0\le t_{+}$. Suppose $t_{+}>0$. Then for every $t\in(0,t_{+})$ there exists a $(K-\lambda_{t})$-bubble, for some $\lambda_{t}$ satisfying (\ref{eq:lambda-bound}). In particular $\lambda_{t}\to\infty$ as $t\to 0^{+}$. Thus the result is proved. If $t_{+}=0$ then $t_{-}<0$ and one argues in a similar way.
\qed

Theorem \ref{T:existence4} holds just for a sequence $|\lambda_{n}|\to\infty$ and not for every large $|\lambda|$ because the set of Lagrange multipliers for constrained minimizers of the isovolumetric problems (\ref{eq:isovolumetric}) in principle could contain gaps.

\section{The isoperimetric problem}
In this section we are mainly interested in the study of the minimization problem defined by (\ref{eq:K-isoperimetric-inequality}).  We assume that $K\in C^{1}(\R^{3})$ satisfies $(K_{2})$--$(K_{3})$ or, as an alternative, $(K_{1})$, $(K_{2})$, and $(K_{4})$. Moreover we suppose that $K\le 0$ on $\R^{3}$. If $K\equiv 0$ then (\ref{eq:K-isoperimetric-inequality}) reduces to the classical isoperimetric inequality (\ref{eq:ii}). In this case extremals exist and are explicitly known, as mentioned in the introduction. Thus we may assume $K\not\equiv 0$. Then by Lemma \ref{L:strict-inequality-SK}, we have that $S_{K}(t)<S_{0}(t)$ for all $t>0$. Hence, by Theorem \ref{T:existence1}, for every $t>0$ the isovolumetric problem defined by (\ref{eq:SKt}) admits a minimizer. We also  point out that, by Lemma \ref{L:SK(t)} (ii) and Lemma \ref{L:SK-estimates} (i),
\begin{gather}
\nonumber%\label{eq:last1}
S_{K}=\inf_{t>0}\tilde{S}_{K}(t)\\
\label{eq:last2}
(1-\tfrac{k_{0}}{2})S\le S_{K}\le S,
\end{gather}
where $S_{K}$ and $\tilde{S}_{K}(t)$ are defined in (\ref{eq:K-isoperimetric-inequality}) and (\ref{eq:StildeK-def}), respectively. 
In the following we study the regularity and the asymptotic behavior of the normalized isovolumetric function $t\mapsto\tilde{S}_{K}(t)$ for $t\in(0,\infty)$. 

\begin{Lemma}
\label{L:Lambda}
The function $t\mapsto\tilde{S}_{K}(t)$ is locally Lipschitz-continuous in $(0,\infty)$. Hence it is differentiable almost everywhere. 
\end{Lemma}

\Proof
Let $t_{1},t_{2}>0$ and take any $u\in M_{1}$. Then
\begin{eqnarray}
\nonumber
\E_{t_{1}}(u)-\E_{t_{2}}(u)&=&\int_{\R^{2}}(Q_{K}(t_{1}u)-Q_{K}(t_{2}u))\cdot u_{x}\wedge u_{y}\\
\label{SK(t)-continuous-1}
&=&\int_{\R^{2}}\left(\int_{t_{2}}^{t_{1}}\frac{\partial}{\partial t}[Q_{K}(tu)]~\!dt\right)\cdot u_{x}\wedge u_{y}.
\end{eqnarray}
Recalling that $Q_{K}(p)=m_{K}(p)p$, one has that
$$
\frac{\partial}{\partial t}[Q_{K}(tp)]=m_{K}(tp)p+(\nabla m_{K}(tp)\cdot p)tp
$$
and, since $K(p)=\mathrm{div}~\!Q_{K}(p)=\nabla m_{K}(p)\cdot p+3m_{K}(p)$, one gets
$$
\frac{\partial}{\partial t}[Q_{K}(tp)]=\frac{1}{t}(K(tp)tp-2Q_{K}(tp)).
$$
Then, taking into account that $|K(p)p|\le k_{0}$ and $|Q_{K}(p)|\le\frac{k_{0}}{2}$ (see Remark \ref{R:K-volume}), one infers that
\begin{equation}
\label{SK(t)-continuous-2}
\left|\int_{t_{2}}^{t_{1}}\frac{\partial}{\partial t}[Q_{K}(tu)]~\!dt\right|\le 2k_{0}\left|\log\frac{t_{2}}{t_{1}}\right|.
\end{equation}
Hence, from (\ref{SK(t)-continuous-1}) and (\ref{SK(t)-continuous-2}) it follows that
\begin{equation}
\label{SK(t)-continuous-3}
\E_{t_{1}}(u)-\E_{t_{2}}(u)\le k_{0}\D(u)\left|\log\frac{t_{2}}{t_{1}}\right|.
\end{equation}
Now take a sequence $(u^{n})\subset M_{1}$ such that $\E_{t_{2}}(u^{n})\to\tilde{S}_{K}(t_{2})$. Since $\tilde{S}_{K}(t_{2})=S_{K_{t_{2}}}(1)$ with $K_{t_{2}}(p)=\mathrm{div}~\!Q_{K}(t_{2}p)$, using Lemma \ref{L:minimizing-PS} (ii), we have that 
\begin{equation}
\label{SK(t)-continuous-4}
\D(u^{n})\le\frac{S}{1-\|Q_{K}\|_{\infty}}+o(1)
\end{equation}
where $o(1)\to 0$ as $n\to\infty$. 
Then, by (\ref{SK(t)-continuous-3})--(\ref{SK(t)-continuous-4}), $\tilde{S}_{K}(t_{1})\le\tilde{S}_{K}(t_{2})+C\left|\log\frac{t_{2}}{t_{1}}\right|$ for some constant $C>0$ independent of $t_{1}$ and $t_{2}$. Exchanging $t_{1}$ with $t_{2}$ we finally obtain
$$
|\tilde{S}_{K}(t_{1})-\tilde{S}_{K}(t_{2})|\le C\left|\log\frac{t_{2}}{t_{1}}\right|
$$
which implies local Lipschitz-continuity in $(0,\infty)$. 
\qed

\begin{Lemma}
\label{L:t=0}
One has that $\tilde{S}_{K}(t)\to S$ as $t\to 0$.
\end{Lemma}

\Proof
Let us recall the following inequality, due to Steffen \cite{Ste76} (see also \cite{CaMuARMA} \textsection 2.3):
\begin{equation}
\label{eq:Kii}
|\Q_{t}(u)|\le\frac{\|K_{t}\|_{\infty}}{S^{3/2}}\D(u)^{3/2}\quad\forall u\in\hat{H}^{1}
\end{equation}
where $K_{t}(p)=\mathrm{div}~\!Q_{K}(tp)=tK(tp)$. If $(u^{n})$ is a minimizing sequence for $\tilde{S}_{K}(t)$ then
\begin{equation}
\label{eq:Ki}
\D(u^{n})\le\frac{\tilde{S}_{K}(t)+o(1)}{1-\|Q_{K}\|_{\infty}}\le\frac{S+o(1)}{1-\|Q_{K}\|_{\infty}}
\end{equation}
where $o(1)\to 0$ as $n\to\infty$. Moreover, by (\ref{eq:Kii}) and (\ref{eq:Ki})
$$
\tilde{S}_{K}(t)\ge\D(u^{n})-\frac{t\|K\|_{\infty}}{S^{3/2}}\D(u^{n})^{3/2}+o(1)\ge S-\frac{t\|K\|_{\infty}}{(1-\|Q_{K}\|_{\infty})^{3/2}}+o(1).
$$
Hence the conclusion follows immediately, using also (\ref{eq:last2}). 
\qed

\begin{Lemma}
\label{L:t=infty}
Let $t_{n}\to\infty$ be such that for every $n$ there exists a minimizer for $\tilde{S}_{K}(t_{n})$. Then $\tilde{S}_{K}(t_{n})\to S$.
\end{Lemma}

\Proof
For every $n$ let $U^{n}\in\hat{H}^{1}$ be a minimizer for $\tilde{S}_{K}(t_{n})$, i.e., $\V(U^{n})=1$ and $\E_{t_{n}}(U^{n})=\tilde{S}_{K}(t_{n})$. In particular
\begin{equation}
\label{eq:nabla2-bdd}
\D(U^{n})\le\frac{\tilde{S}_{K}(t_{n})}{1-\|Q_{K}\|_{\infty}}\le\frac{S}{1-\|Q_{K}\|_{\infty}}.
\end{equation}
Since 
$$
1=\V(U^{n})\le\frac{1}{6}\|U^{n}\|_{\infty}\D(U^{n}),
$$ one has that 
$$
\|U^{n}\|_{\infty}\ge\frac{6S}{1-\|Q_{K}\|_{\infty}}=:\delta_{0}>0.
$$
Recall that $U^{n}$ is a $(K_{n}-\lambda_{n})$-bubble, with $K_{n}(p)=\mathrm{div}~\!Q_{K}(t_{n}p)$ and $\lambda_{n}\in\R$. In particular $U^{n}$ is bounded and regular as a map on $\S^{2}$ and there exists $U^{n}(\infty)=\lim_{|z|\to\infty}U^{n}(z)$. By the conformal invariance, without changing notation, we may assume that $|U^{n}(\infty)|=\|U^{n}\|_{\infty}$ and that $\{|U^{n}|<\delta_{0}/2\}\subset\disc$ where in general
$$
\{|U^{n}|<\delta\}:=\{z\in\R^{2}~|~|U^{n}(z)|<\delta\}
$$
and $\disc$ denotes the open unit disc. For every $\delta\in(0,\delta_{0})$ let
\begin{equation}
\label{eq:notation1}
A_{n}(\delta):=\int_{\{|U^{n}|<\delta\}}|\nabla U^{n}|^{2}\quad(n\in\N)\quad\text{and}\quad A(\delta):=\liminf_{n\to\infty}A_{n}(\delta).
\end{equation}
The following technical result holds. 
\begin{Lemma}
\label{L:A(delta)}
One has that $A(\delta)\to 0$ as $\delta\to 0^{+}$.
\end{Lemma}

\noindent
Let us complete the proof of Lemma \ref{L:t=infty}. Since $\tilde{S}_{K}(t)\le S$ for every $t>0$, it is enough to show that $\liminf\tilde{S}_{K}(t_{n})\ge S$. Fix $\varepsilon>0$. 
By Lemma \ref{L:A(delta)}, there exists $\delta_{\varepsilon}>0$ such that $A(\delta_{\varepsilon})\le\varepsilon$, namely,
\begin{equation}
\label{eq:limif}
\liminf_{n\to\infty}\int_{\{|U^{n}|<\delta_{\varepsilon}\}}|\nabla U^{n}|^{2}\le\varepsilon.
\end{equation}
Then we estimate
$$
|\Q_{t_{n}}(U^{n})|\le\int_{\{|U^{n}|\ge\delta_{\varepsilon}\}}|Q_{K}(t_{n}U^{n})|~\!|U^{n}_{x}\wedge U^{n}_{y}|+\frac{\|Q_{K}\|_{\infty}}{2}\int_{\{|U^{n}|<\delta_{\varepsilon}\}}|\nabla U^{n}|^{2}.
$$
Using (\ref{eq:Q(infty)=0}), (\ref{eq:nabla2-bdd}), and (\ref{eq:limif}), we obtain that
$$
\liminf_{n\to\infty}|\Q_{t_{n}}(U^{n})|\le\frac{\|Q_{K}\|_{\infty}\varepsilon}{2}.
$$
Thus, by the arbitrariness of $\varepsilon>0$, one deduces that $\Q_{t_{n}}(U^{n})\to 0$ for a subsequence. Consequently, since $\V(U^{n})=1$,  
$$
\tilde{S}_{K}(t_{n})=\D(U^{n})+\Q_{t_{n}}(U^{n})\ge S-|\Q_{t_{n}}(U^{n})|=S+o(1)\quad\text{as $n\to\infty$}
$$
and we are done.
\qed

\noindent
\emph{Proof of Lemma \ref{L:A(delta)}.} We know that $U^{n}$ is a minimizer for $\tilde{S}_{K}(t_{n})$ and solves
\begin{equation}
\label{eq:Un-system}
-\Delta U^{n}+t_{n}K(t_{n}U^{n})U^{n}_{x}\wedge U^{n}_{y}=\lambda_{n}U^{n}_{x}\wedge U^{n}_{y}\quad\text{on }\R^{2}
\end{equation}
for some $\lambda_{n}$. Then $t_{n}U^{n}$ is a minimizer for $S_{K}(t_{n}^{3})$ and is a $\big(K-\frac{\lambda_{n}}{t_{n}}\big)$-bubble. By (\ref{eq:lambda-bound}) we infer that
\begin{equation}
\label{eq:lambdan-bound}
c_{1}:=\frac{(2-k_{0})^{2}S}{3(2+k_{0})}\le\lambda_{n}\le\frac{2(2+k_{0})S}{3(2-k_{0})}=:c_{2}.
\end{equation}
For every $\delta\in(0,\delta_{0}/2)$ let 
$$
\phi_{\delta}(s):=\left\{\begin{array}{ll}1&\text{as $0\le s\le\delta$}\\ \frac{2\delta}{s}-1&\text{as $\delta< s\le2\delta$}\\ 0&\text{as $s>2\delta.$}\end{array}\right.
$$
Set $u^{n}:=\phi_{\delta}(|U^{n}|)U^{n}$. Since $\phi_{\delta}$ is Lipschitz-continuous and $\{|U^{n}|<\delta_{0}/2\}\subset\disc$, one has that $u^{n}\in H^{1}_{0}(\disc,\R^{3})$. We test (\ref{eq:Un-system}) with $u^{n}$ and we find
\begin{equation}
\label{eq:B0}
\begin{split}
\int_{\disc}\nabla U^{n}\cdot\nabla u^{n}+&\int_{\disc}\phi_{\delta}(|U^{n}|)K(t_{n}U^{n})t_{n}U^{n}\cdot U^{n}_{x}\wedge U^{n}_{y}\\
&=\lambda_{n}\int_{\disc}\phi_{\delta}(|U^{n}|)U^{n}\cdot U^{n}_{x}\wedge U^{n}_{y}.
\end{split}
\end{equation}
Let us estimate each term in (\ref{eq:B0}) as follows. Firstly
$$
\int_{\disc}\nabla U^{n}\cdot\nabla u^{n}=\int_{\{|U^{n}|<\delta\}}|\nabla U^{n}|^{2}+\int_{\{\delta<|U^{n}|<2\delta\}}\nabla U^{n}\cdot\nabla u^{n}.
$$
On the set $\{\delta<|U^{n}|<2\delta\}$, with direct computations, one finds that
$$
\nabla U^{n}\cdot\nabla u^{n}=-2\delta~\!\frac{(U^{n}\cdot U^{n}_{x})^{2}+(U^{n}\cdot U^{n}_{y})^{2}}{|U^{n}|^{3}}+2\delta~\!\frac{|\nabla U^{n}|^{2}}{|U^{n}|}-|\nabla U^{n}|^{2}\ge -|\nabla U^{n}|^{2}
$$
and then
\begin{equation}
\label{eq:B1}
\int_{\disc}\nabla U^{n}\cdot\nabla u^{n}\ge\int_{\{|U^{n}|<\delta\}}|\nabla U^{n}|^{2}-\int_{\{\delta<|U^{n}|<2\delta\}}\nabla U^{n}\cdot\nabla u^{n}=2A_{n}(\delta)-A_{n}(2\delta),
\end{equation}
according to the notation introduced in (\ref{eq:notation1}). Secondly, by $(K_{1})$,
\begin{equation}
\label{eq:B2}
\left|\int_{\disc}\phi_{\delta}(|U^{n}|)K(t_{n}U^{n})t_{n}U^{n}\cdot U^{n}_{x}\wedge U^{n}_{y}\right|\le k_{0}\int_{\{|U^{n}|<2\delta\}}|U^{n}_{x}\wedge U^{n}_{y}|\le\frac{k_{0}}{2}~\!A_{n}(2\delta).
\end{equation}
In addition
\begin{equation}
\label{eq:B3}
\left|\int_{\disc}\phi_{\delta}(|U^{n}|)U^{n}\cdot U^{n}_{x}\wedge U^{n}_{y}\right|\le \int_{\{|U^{n}|<2\delta\}}|U^{n}|~\!|U^{n}_{x}\wedge U^{n}_{y}|\le\delta~\!A_{n}(2\delta).
\end{equation}
Hence, from (\ref{eq:B0})--(\ref{eq:B3}) it follows that
$$
2A_{n}(\delta)\le\left(1+\frac{k_{0}}{2}+\lambda_{n}\delta\right)A_{n}(2\delta)\quad\forall n\in\N
$$
and then, using also (\ref{eq:lambdan-bound}), 
\begin{equation}
\label{eq:B4}
2A(\delta)\le\left(1+\frac{k_{0}}{2}+\overline{\lambda}\delta\right)A(2\delta)
\end{equation}
for some $\overline\lambda>0$ and for every $\delta\in(0,\delta_{0}/2)$. The mappings $\delta\mapsto A_{n}(\delta)$ are non-negative and non-decreasing, and the same holds for the mapping $\delta\mapsto A(\delta)$. In particular there exists 
$$
\lim_{\delta\to 0^{+}}A(\delta)=:A\ge 0
$$
and, by (\ref{eq:B4}), $2A\le\left(1+\frac{k_{0}}{2}\right)A$. Since $k_{0}<2$ we conclude that $A=0$.
\qed

\begin{Theorem}
\label{T:Lambda}
For every $t\in(0,t_{+})$ let $\Lambda(t)$ be the set of Lagrange multipliers for minimizers of $S_{K}(t)$, i.e., $\Lambda(t)=\{\lambda\in\R~|~\exists~(K-\lambda)\text{-bubble }U\in M_{t}\text{ such that }\E(U)=S_{K}(t)\}$. Then:
\begin{itemize}[leftmargin=20pt]
\item[(i)]
for every $t\in(0,t_{+})$ the set $\Lambda(t)$ is compact, $\Lambda(t)\subset[c_{1}t^{-1/3},c_{2}t^{-1/3}]$ with $c_{1}$ and $c_{2}$ defined in (\ref{eq:lambdan-bound}). Moreover  
\begin{equation}
\label{eq:min-max-Lambda}
\begin{split}
&\limsup_{\varepsilon\to 0^{+}}\frac{S_{K}(t+\varepsilon)-S_{K}(t)}{\varepsilon}\le\min\Lambda(t),\\
&\liminf_{\varepsilon\to 0^{-}}\frac{S_{K}(t+\varepsilon)-S_{K}(t)}{\varepsilon}\ge\max\Lambda(t).
\end{split}
\end{equation}
\item[(ii)]
For a.e. $t\in(0,t_{+})$ there exists the derivative $S_{K}'(t)$ and  $\Lambda(t)=\{S'_{K}(t)\}$.
\end{itemize}
\end{Theorem}

\Proof
(i) By (\ref{eq:lambda-bound}), $\Lambda(t)$ is bounded. Let $(\lambda_{n})\subset\Lambda(t)$ be such that $\lambda_{n}\to\lambda$, $\lambda_{n}\ne\lambda$. Then there is a sequence $(U^{n})\subset\hat{H}^{1}$ of minimizers for $S_{K}(t)$ and each $U^{n}$ is a $(K-\lambda_{n})$-bubble. Since $t\in(0,t_{+})$, we have that $S_{K}(\tau)<S_{0}(\tau)$ for every $\tau\in(0,t]$ (Lemma \ref{L:strict-inequality-SK}). The sequence $(U^{n})$ satisfies (\ref{eq:almost-solution}). Indeed
\begin{equation*}\begin{split}
|\E'(U^{n})[\varphi]-\lambda\V'(U^{n})[\varphi]|&=|\lambda_{n}-\lambda|~\!|\V'(U^{n})[\varphi]|\\
&\le C|\lambda_{n}-\lambda|~\!\|\nabla U^{n}\|_{2}^{2}\|\nabla\varphi\|_{2}\quad\forall\varphi\in C^{\infty}_{c}(\R^{2},\R^{3})
\end{split}\end{equation*}
with $C>0$ independent of $n$ (see Remark \ref{R:volume}). Moreover $c_{1}\le\|\nabla U^{n}\|_{2}\le c_{2}$ for some constants $0<c_{1}\le c_{2}<\infty$ (Lemma \ref{L:minimizing-PS}), and $C^{\infty}_{c}(\R^{2},\R^{3})$ is dense in $\hat{H}^{1}$. Hence $\|\E'(U^{n})-\lambda\V'(U^{n})\|_{\hat{H}^{-1}}\to 0$. Therefore we can apply Lemma \ref{L:PS-decomposition} and repeating the proof of Lemma \ref{L:minimizer}, we infer that the decomposition of $(U^{n})$ according to Lemma \ref{L:PS-decomposition} in fact is made by just one $(K-\lambda)$-bubble $U\in\hat{H}^{1}$. Moreover $\E(U^{n})\to\E(U)$ and $\V(U^{n})\to\V(U)$. In particular $U$ is a minimizer for $S_{K}(t)$. Thus we proved that $\lambda\in\Lambda(t)$, namely $\Lambda(t)$ is closed. Now let us prove (\ref{eq:min-max-Lambda}).
For every $\lambda\in\Lambda(t)$ there exists a $(K-\lambda)$-bubble $U\in\hat{H}^{1}$ which is a minimizer for $S_{K}(t)$. Let $u^{\varepsilon}=\sqrt[3]{1+\frac{\varepsilon}{t}}~\!U$. Then $\V(u^{\varepsilon})=t+\varepsilon$, $S_{K}(t+\varepsilon)\le\E(u^{\varepsilon})$, and 
\begin{equation*}\begin{split}
\limsup_{\varepsilon\to 0^{+}}\frac{S_{K}(t+\varepsilon)-S_{K}(t)}{\varepsilon}&\le\lim_{\varepsilon\to 0}\frac{\E(u^{\varepsilon})-\E(U)}{\varepsilon}\\
&=\lim_{\varepsilon\to 0}\frac{\sqrt[3]{1+\frac{\varepsilon}{t}}-1}{\varepsilon}\lim_{s\to 1}\frac{\E(sU)-\E(U)}{s-1}=\frac{\E'(U)[U]}{3t}
\end{split}\end{equation*}
Since $U$ is a $(K-\lambda)$-bubble, we have that $\E'(U)[U]=3\lambda t$ and thus we get
$$
\limsup_{\varepsilon\to 0^{+}}\frac{S_{K}(t+\varepsilon)-S_{K}(t)}{\varepsilon}\le\lambda.
$$
In a similar way we can show the opposite inequality for the $\liminf$ as $\varepsilon\to 0^{-}$. Thus (\ref{eq:min-max-Lambda}) holds.
\medskip

\noindent
(ii)
By Lemma \ref{L:Lambda}, the isovolumetric mapping $t\mapsto S_{K}(t)=t^{2/3}\tilde{S}_{K}(t^{1/3})$ is differentiable a.e. Hence, if there exists the derivative $S'_{K}(t)$, by (\ref{eq:min-max-Lambda}) the set $\Lambda(t)$ is a singleton and its unique element is $S'_{K}(t).\quad\square$
\medskip

\noindent
\emph{Proof of Theorem \ref{T:existence3}.}
The normalized isovolumetric function $t\mapsto\tilde{S}_{K}(t)$ is continuous in $(0,\infty)$, as stated in Lemma \ref{L:Lambda}. Then, by (\ref{eq:last2}) and by Lemmas \ref{L:t=0} and \ref{L:t=infty}, there exists $t_{0}>0$ such that $\tilde{S}_{K}(t_{0}^{1/3})=\inf_{t>0}\tilde{S}_{K}(t)=S_{K}$. Let $U\in\hat{H}^{1}$ be a minimizer for the isovolumetric problem (\ref{eq:SKt}) with $t=t_{0}$. Then one easily checks that 
$$
\frac{\E(U)}{\V(U)^{2/3}}=S_{K}
$$
namely $U$ is a minimizer for (\ref{eq:K-isoperimetric-inequality}). In particular $U$ is a $(K-\lambda)$-bubble for some $\lambda>0$. Since $t^{2}\tilde{S}_{K}(t)=S_{K}(t^{3})$ (Lemma \ref{L:SK(t)} (ii)), we have that
%\begin{equation}
\begin{eqnarray}
\nonumber
&&\limsup_{\varepsilon\to 0^{+}}\frac{\tilde{S}_{K}(t_{0}^{1/3}+\varepsilon)-\tilde{S}_{K}(t_{0}^{1/3})}{\varepsilon}+\frac{2\tilde{S}_{K}(t_{0}^{1/3})}{t_{0}^{1/3}}
\\
\label{eq:limsup}
&&\qquad\qquad\le 3\displaystyle\limsup_{\delta\to 0^{+}}\frac{{S}_{K}(t_{0}+\delta)-{S}_{K}(t_{0})}{\delta}\le 3\min\Lambda(t_{0})
\end{eqnarray}
%\end{equation}
in view of (\ref{eq:min-max-Lambda}). Similarly one can obtain
\begin{equation}
\label{eq:liminf}
\liminf_{\varepsilon\to 0^{-}}\frac{\tilde{S}_{K}(t_{0}^{1/3}+\varepsilon)-\tilde{S}_{K}(t_{0}^{1/3})}{\varepsilon}+\frac{2\tilde{S}_{K}(t_{0}^{1/3})}{t_{0}^{1/3}}\ge 3\max\Lambda(t_{0}).
\end{equation}
Since $\tilde{S}_{K}(t_{0}^{1/3})=\min_{t>0}\tilde{S}_{K}(t)=S_{K}$, from (\ref{eq:limsup}) and (\ref{eq:liminf}) it follows that $\min\Lambda(t_{0})\ge\frac{2}{3}S_{K}t_{0}^{-1/3}$ and $\max\Lambda(t_{0})\le\frac{2}{3}S_{K}t_{0}^{-1/3}$, respectively. Thus $\Lambda(t_{0})=\big\{\tfrac{2}{3}S_{K}t_{0}^{-1/3}\big\}$, that is  $\lambda=\frac{2}{3}S_{K}\V(U)^{-1/3}$. 
\qed

\section{A nonexistence result for isovolumetric problems}

Some tools introduced in the previous section can be also used to show a nonexistence result for the isovolumetric problems (\ref{eq:isovolumetric}). Such a result has a counterpart in the context of the $H$-bubble problem (see \cite{CaMuARMA}, \textsection 6). Here we assume the strict inequality $K<0$ on $\R^{3}$.

\begin{Theorem}
\label{T:nonexistence}
Let $K\in C^{0}(\R^{3})$ satisfy $(K_{1})$--$(K_{2})$. If $K<0$ on $\R^{3}$, then there exists $\varepsilon>0$ such that $\tilde{S}_{K}(t)=S$ and the isovolumetric problem (\ref{eq:isovolumetric}) has no minimizer for all $t\in(-\varepsilon,0)$.
\end{Theorem}

\Proof
Firstly let us prove that no minimizer for $S_{K}(t)$ exists as $t<0$ with small $|t|$. Arguing by contradiction, assume that in correspondence of a sequence $t_{n}\to 0^{-}$, for every $n$ there exists a minimizer $U^{n}\in\hat{H}^{1}$ for $S_{K}(t_{n})$. Setting $\tau_{n}=\sqrt[3]{t_{n}}$ and $u^{n}=\frac{1}{\tau_{n}}U^{n}$, each $u^{n}$ turns out to be a minimizer for $\tilde{S}_{K}(\tau_{n})$. Moreover, since $\tau_{n}\to 0$, $(u^{n})$ is a minimizing sequence for the isoperimetric problem defined by
$$
S=\inf\{\D(u)~|~u\in\hat{H}^{1},~\V(u)=1\}
$$
because $\V(u^{n})=1$ and by an application of (\ref{eq:Kii}) and (\ref{eq:Ki}). By known results (see Lemma 2.1 in \cite{CaMuARMA}), there exist a sequence of conformal mappings $g_{n}\colon\S^{2}\to\S^{2}$ and a sequence $(p_{n})\subset\R^{3}$ such that $u^{n}\circ g_{n}-p_{n}\to-\omega$ strongly in $\hat{H}^{1}$, where $\omega$ is the standard parametrization of the unit sphere, defined in Remark \ref{R:sphere}. In fact, by conformal invariance, the function $\tilde{U}^{n}=u^{n}\circ g_{n}$ is also a minimizer for $\tilde{S}_{K}(\tau_{n})$, hence is a $(\tilde{K}_{n}-\lambda_{n})$-bubble, where $\tilde{K}_{n}(p)=\tau_{n}K(\tau_{n}p)$ and $\lambda_{n}$ is bounded. Using an $\varepsilon$-regularity argument (see, e.g., the last part of the proof of Theorem 6.3 in \cite{CaMuARMA} and the references therein), we can show that $\tilde{U}^{n}-p_{n}\to-\omega$ in $C^{1}(\S^{2},\R^{3})$. This implies that for $n$ large enough, $\tilde{U}^{n}$ is an embedded parametric surface, bounding a domain $A_{n}\subset\R^{3}$ and 
$$
\Q_{t_{n}}(\tilde{U}^{n})=\int_{A_{n}}\tilde{K}_{n}(p)~\!dp.
$$
Since $K<0$ on $\R^{3}$ and $\tau_{n}<0$, we obtain that 
$$
\tilde{S}_{K}(\tau_{n})=\D(\tilde{U}^{n})+\Q_{t_{n}}(\tilde{U}^{n})>\D(\tilde{U}^{n})\ge S,
$$
namely $S_{K}(t_{n})>St_{n}^{2/3}$, contrary to Lemma \ref{L:SK-estimates} (i). Now we show that $\tilde{S}_{K}(t)=S$ as $t<0$ with small $|t|$. Again we argue by contradiction, assuming that $\tilde{S}_{K}(\tau_{n})<S$ along a sequence $\tau_{n}\to 0^{-}$. Then $S_{K}(t_{n})<St_{n}^{2/3}$ for every $n$, where $t_{n}=\tau_{n}^{3}$. Reasoning as in the proof of Lemma \ref{L:minimizer}, for fixed $n$, we find a decomposition of $t_{n}=\sum_{i\in I_{n}}t_{n,i}+\sum_{j\in J_{n}}t_{n,j}$ with $t_{n,i},t_{n,j}\in[t_{n},0)$, $I_{n}$ and $J_{n}$ finite sets of indices, with $S_{K}(t_{n,i})$ admitting a minimizer, and 
$$
S_{K}(t_{n})=\sum_{i\in I_{n}}S_{K}(t_{n,i})+\sum_{j\in J_{n}}S_{0}(t_{n,j}).
$$
Notice that the assumptions $(K_{1})$ and $(K_{2})$ are enough for this part of the argument. We claim that $I_{n}\ne\varnothing$. If not, then we reach a contradiction because
$$
St_{n}^{2/3}>S_{K}(t_{n})=\sum_{j\in J_{n}}S_{0}(t_{n,j})=\sum_{j\in J_{n}}St_{n,j}^{2/3}\ge S\bigg(\sum_{j\in J_{n}}t_{n,j}\bigg)^{2/3}=St_{n}^{2/3}.
$$
Thus we proved that there exists $t'_{n}\in(t_{n},0)$ such that the isovolumetric problem defined by $S_{K}(t'_{n})$ admits a minimizer. Then we apply the first part of the proof to reach a contradiction. 
\qed
\bigskip

\noindent
\textbf{Acknowledgements.} Work partially supported by the PRIN-2012-74FYK7 Grant ``Variational and perturbative aspects of nonlinear differential problems'', by the project ERC Advanced Grant 2013 n. 339958 ``Complex Patterns for Strongly Interacting Dynamical Systems - COMPAT'', and by the Gruppo Nazionale per l'Analisi Mate\-ma\-tica, la Probabilit\`a e le loro Applicazioni (GNAMPA) of the Istituto Nazionale di Alta Mate\-ma\-tica (INdAM).

\end{document}